\pdfoutput=1
\RequirePackage{ifpdf}
\ifpdf 
\documentclass[pdftex]{sigma}
\else
\documentclass{sigma}
\fi

\numberwithin{equation}{section}

\newtheorem{Theorem}{Theorem}[section]
\newtheorem*{Theorem*}{Theorem}
\newtheorem{Corollary}[Theorem]{Corollary}
\newtheorem{Lemma}[Theorem]{Lemma}
\newtheorem{Proposition}[Theorem]{Proposition}
 { \theoremstyle{definition}
\newtheorem{Definition}[Theorem]{Definition}

\newtheorem{Remark}[Theorem]{Remark} }

\newcommand{\E}[1]{\widetilde{e}_{#1}}

\newcommand{\boldcross}{\includegraphics{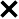}}

\def\a{\alpha}

\def\C{\mathbb{C}}

\def\e{\epsilon}
\def\E{\mathbb{E}}

\def\i{\infty}

\def\l{\lambda}

\def\o{\omega}

\def\P{\mathbb{P}}
\def\Pe{\mathcal{P}}
\def\1{\bf{1}}
\def\R{\mathbb{R}}

\def\S{\mathcal{S}}

\def\t{\tau}

\def\Z{\mathbb{Z}}

\begin{document}


\newcommand{\arXivNumber}{2106.11913}

\renewcommand{\PaperNumber}{064}

\FirstPageHeading

\ShortArticleName{Identity between Restricted Cauchy Sums}

\ArticleName{Identity between Restricted Cauchy Sums\\ for the $\boldsymbol{q}$-Whittaker and Skew Schur Polynomials}

\Author{Takashi IMAMURA~$^{\rm a}$, Matteo MUCCICONI~$^{\rm b}$ and Tomohiro SASAMOTO~$^{\rm c}$}

\AuthorNameForHeading{T.~Imamura, M.~Mucciconi and T.~Sasamoto}

\Address{$^{\rm a)}$~Department of Mathematics and Informatics, Chiba University, Chiba, 263-8522 Japan}
\EmailD{\href{mailto:imamura@math.s.chiba-u.ac.jp}{imamura@math.s.chiba-u.ac.jp}}

\Address{$^{\rm b)}$~Department of Mathematics, University of Warwick, Coventry, CV4 7HP, UK}
\EmailD{\href{mailto:matteomucciconi@gmail.com}{matteomucciconi@gmail.com}}

\Address{$^{\rm c)}$~Department of Physics, Tokyo Institute of Technology, Tokyo, 152-8551 Japan}
\EmailD{\href{mailto:sasamoto@phys.titech.ac.jp}{sasamoto@phys.titech.ac.jp}}

\ArticleDates{Received December 20, 2023, in final form July 02, 2024; Published online July 16, 2024}

\Abstract{The Cauchy identities play an important role in the theory of symmetric functions. It is known that Cauchy sums for the $q$-Whittaker and the skew Schur polynomials produce the same factorized expressions modulo a $q$-Pochhammer symbol. We consider the sums with restrictions on the length of the first rows for labels of both polynomials and prove an identity which relates them. The proof is based on techniques from integrable probability: we rewrite the identity in terms of two probability measures: the $q$-Whittaker measure and the periodic Schur measure. The relation follows by comparing their Fredholm determinant formulas.}


\Classification{05A19; 05E05; 60J10}

\section{Introduction}
\label{s1}

\subsection{The identity and connections to integrable probability}
 Symmetric polynomials have been of much interest in various fields of mathematics such as representation theory, combinatorics etc. Among them the Schur polynomials $s_{\l}(a_1,\dots,a_N)$ occupy an important place because of their simplicity, depth and wide applications. One of their most well-known properties
 is the Cauchy identity,
\begin{align}
\label{1-1}
 \sum_{\l\in\Pe}s_{\l}(a_1,\dots,a_N)s_{\l}(b_1,\dots,b_M)
 =
 \prod_{i=1}^N\prod_{j=1}^M\frac{1}{1-a_ib_j},
\end{align}
 where $\Pe$ is the set of partitions.
 Along with a number of related summation identities, \eqref{1-1} provides important insights into various field of mathematics such as combinatorics \cite{Sagan2000symmetricgroup}, representation theory \cite{Bump2013Liegroups}, integrable probability \cite{BorodinGorin2016IntegrableProbability}, etc.

 The Cauchy identity for the Schur function \eqref{1-1} has been generalized in various different ways. In this paper, we focus on two generalizations of \eqref{1-1} in which the Schur polynomials in the left-hand side are replaced by the $q$-Whittaker functions $P_{\l}(a_1,\dots,a_N)$, $Q_{\l}(b_1,\dots,b_M)$ and the skew Schur polynomials $s_{\l/\rho}(a_1,\dots,a_N)$. These are
 \begin{align}
 \label{1-3}
 &\sum_{\mu\in\Pe}P_{\mu}(a_1,\dots,a_N)
 Q_{\mu}(b_1,\dots,b_M)
 =
 \prod_{i=1}^N\prod_{j=1}^M\frac{1}{(a_ib_j;q)_{\i}},
 \\
 \label{1-2}
 &\sum_{\substack{\l,\rho\in\Pe\\
 \rho\subset \l
 }}q^{|\rho|}s_{\l/\rho}(a_1,\dots,a_N)s_{\l/\rho}(b_1,\dots,b_M)
 =
 \frac{1}{(q;q)_\i}
 \prod_{i=1}^N\prod_{j=1}^M\frac{1}{(a_ib_j;q)_{\i}},
 \end{align}
 where $|\rho|:=\rho_1+\rho_2+\cdots$, $\rho\subset\l$ means $\rho_i\le \l_i$ for any $i=1,2,\dots$ and $(x;q)_{\i}$ is defined by~\eqref{a2-1}. The $q$-Whittaker functions $P_{\l}(a_1,\dots,a_N)$ are a special case of the Macdonald symmetric polynomials $P_{\l}(a;q;t)$ defined in \cite[Chapter~VI]{macdonald1998symmetric}, in which out of two parameters $(q,t)$ of the latter one puts $t=0$. The dual polynomial $Q_{\mu}(b_1,\dots,b_M)$ is defined as
 $
 \prod_{j=1}^M (q;q)_{\mu_j-\mu_{j+1}}^{-1}\cdot P_{\mu}(b_1,\dots,b_M)$,
 where $(x;q)_n$ is given by \eqref{a2-1d} and we set $\mu_{M+1}\equiv 0$. When $q=0$, both $P_{\l}(a_1,\dots,a_N)$ and~${Q_{\l}(a_1,\dots,a_N)}$ reduce to the Schur polynomial and \eqref{1-2} becomes \eqref{1-1}. The skew Schur polynomials \smash{$s_{\l/\rho}(a_1,\dots,a_N)$} are labeled by skew partitions $\l/\rho$. For a precise definition, see, e.g., \cite[Chapter~I]{macdonald1998symmetric}. When $q=0$, since we have $\rho=\varnothing$ due to the weight $q^{|\rho|}$, they reduce to (non-skew) Schur polynomial $s_{\l}(a_1,\dots,a_N)$ and \eqref{1-3} becomes \eqref{1-1}. Proofs of \eqref{1-3} and~\eqref{1-2} can be found in \cite[Sections~VI.4 and~I.5.28\,(a)]{macdonald1998symmetric}.

 We immediately see that the right-hand sides of \eqref{1-3} and \eqref{1-2} are exactly the same up to a factor $1/(q;q)_{\i}$. Combining this fact with \eqref{a2-8}, which rewrites $1/(q;q)_{\i}$ as a sum over partitions $\nu$, we have
 \begin{align}
 \label{1-4}
 \sum_{\mu,\nu\in\Pe}q^{|\nu|}P_{\mu}(a_1,\dots,a_N)
 Q_{\mu}(b_1,\dots,b_M)
 =
 \sum_{\substack{\l,\rho\in\Pe\\ \rho\subset\l}}q^{|\rho|}s_{\l/\rho}(a_1,\dots,a_N)s_{\l/\rho}(b_1,\dots,b_M).
 \end{align}
 This identity suggests that there may exist deeper connections between the $q$-Whittaker functions and the skew Schur functions resulting in possible refinements of \eqref{1-4}.
 Indeed, in this paper, we will prove such an identity as stated in the following theorem.

 \begin{Theorem} \label{t1}
 For any $n\in\Z_{\ge 0}$, we have
 \begin{gather}
 \sum_{\stackrel{\mu,\nu\in\Pe}{\mu_1+\nu_1\le n}}q^{|\nu|}P_{\mu}(a_1,\dots,a_N)
 Q_{\mu}(b_1,\dots,b_M)\nonumber\\
 \qquad=
 \sum_{\stackrel{\l,\rho\in\Pe}{\rho\subset\l,\,\l_1\le n}}
 q^{|\rho|}s_{\l/\rho}(a_1,\dots,a_N)s_{\l/\rho}(b_1,\dots,b_M). \label{t1-1}
 \end{gather}
 \end{Theorem}
 Summations of this type with restrictions on the length of the first row (mostly for the Schur case) have attracted attention in the theory of symmetric functions in connection with combinatorics and representation theory, see for instance \cite{Baik-Rains2001algebraic,Gessel1990Cauchysum,Rains1998IncreasingSubsequence,Rains-Warnaar2015boundedLittlewoodidentities}.
 In this paper, the proof of Theorem~\ref{t1} will be based on connections of quantities on both sides to integrable probability. In our next paper \cite{ImamuraMucciconiSasamoto2021skewRSK}, we will develop a combinatorial theory which allows a bijective proof of~\eqref{t1-1} and related identities.

 \begin{Remark} \label{rem:intro}
 In Theorem~\ref{t1}, our main identity is stated for complex variables $a_1,\dots,a_N$, $b_1,\dots,b_M$, but one easily sees that the same holds also when they are arbitrary specialization of symmetric functions. In order to establish \eqref{t1-1}, it suffices to prove it for $a_i$, $b_j$ varying in an continuous interval, which is in fact what we will prove below using integrable probabilistic arguments.
 \end{Remark}

 Let us rewrite \eqref{t1-1} in a probabilistic language. The probability measures on $\Pe$ related to~\eqref{1-3} and \eqref{1-2} are of much interest in integrable probability.
 They are defined as
 \begin{gather}
 M_{q\text{W}}(\mu)=P_{\mu}(a_1,\dots,a_N)Q_{\mu}(b_1,\dots,b_M)/Z_{q\text{W}}, \label{s13}
 \\
 M_{\text{pS}}(\l)=
 \sum_{\substack{\rho\in\Pe\\ \rho \subset\l }}q^{|\rho|}s_{\l/\rho}(a_1,\dots,a_N)s_{\l/\rho}(b_1,\dots,b_M)/Z_{\text{pS}}, \label{1-8}
 \end{gather}
 where we assumed $q \in (0,1)$ $a_i,b_j>0$, with $a_ib_j<1$ for $i=1,\dots,N$ and $j=1,\dots,M$ and the normalization constants $Z_{q\text{W}}$ and $Z_{\text{pS}}$ are given by the right-hand sides of the Cauchy identities~\eqref{1-3} and \eqref{1-2} respectively.
 $M_{q\text{W}}(\mu)$ is called the \emph{$q$-Whittaker measure}
 while $M_{\text{pS}}(\l)$ is called the \emph{periodic Schur measure}.

 The $q$-Whittaker measure was introduced in \cite{BorodinCorwin2014Mac}. It has been playing an important role in the development of integrable probability in the past decade. Its marginals are known to describe one point functions of certain particle systems discretizing the KPZ equation. On the other hand, the periodic Schur measure was introduced in \cite{Borodin2007PeriodicSchur}. It can be regarded as a model of lozenge tilings in a cylindric domain and free fermions at positive temperature \cite{BeteaBouttier2019PeriodicSchur}.

 We also introduce two random variables $\chi$ and $S$, which are independent of other random variables and respectively have distribution,
 \begin{align}
\label{t11-3}
\P(\chi=n)=\frac{q^n}{(q;q)_n}(q;q)_{\i},
\qquad
n=0,1,2,\dots
\end{align}
and for a fixed parameter $t \in \mathbb{R}_{>0}$,
 \begin{align}
 \label{1-11}
 \P(S=\ell):=\frac{t^{\ell} q^{\ell^2/2}}{(q;q)_\i
 \theta \bigl(-tq^{1/2}\bigr)},
 \qquad
\ell \in \mathbb{Z}.
\end{align}
Here the definitions of the $q$-Pochhammer symbols
$(x;q)_\i$ and $(x;q)_n$ are given by \eqref{a2-1} and~\eqref{a2-1d} respectively whereas the function $\theta(x)$ is defined below \eqref{a2-6}. Being proportional to $q^{\ell^2/2},~\ell\in\mathbb{Z}$, the distribution defined by \eqref{1-11} is sometimes called a discrete Gaussian distribution.

The random variable $\chi$ is related to the partition $\nu$ in \eqref{t1-1}.
Namely, thanks to the summation identity \eqref{a2-7}, if a partition $\nu$ is sampled with probability proportional to $q^{|\nu|}$, then its first row~$\nu_1$ follows the same law as $\chi$. On the other hand the random variable $S$ is introduced as it unveils a~determinantal structure in the periodic Schur process. This was observed by Borodin in \cite{Borodin2007PeriodicSchur}, who noticed that shifting by $S$ all entries of a~random partition $\lambda$ sampled with $M_{\mathrm{pS}}$ law produces a~determinantal point process called \emph{shift-mixed periodic Schur measure}.

With these preparations, Theorem~\ref{t1} can be restated as follows.
 \begin{Theorem} \label{thm:intro_prob}\label{t12}
 Let $\mu_1$, $\lambda_1$ be the first rows of partitions $\mu$, $\lambda$, distributed respectively according to the $q$-Whittaker measure $M_{q{\rm W}}(\mu)$ \eqref{s13} and to the periodic Schur measure $M_{{\rm pS}}(\l)$ \eqref{1-8}. Let also $\chi$ and $S$ be independent random variables distributed as \eqref{t11-3} and \eqref{1-11}, respectively.
 The following equalities hold and they are equivalent. For $n\in\Z_{\ge 0}$,
 we have
 \begin{align}
 \label{t12-1}
& {\rm (a)}\quad \P (\mu_1+\chi\le n )=
 \P (\l_1\le n ),
\\
\label{t12-3} & {\rm (b)} \quad
 \E\left[\frac{1}{\bigl(-t q^{\frac12+n-\mu_1};q\bigr)_{\i}}\right]
 =
 \P(\l_1+S\le n).
\end{align}
\end{Theorem}
 Dividing both sides of \eqref{t1-1} by the right-hand side by $Z_{\textrm{pS}}$ and noting $Z_{\textrm{pS}}=Z_{q\textrm{W}}Z_{\chi}$, where~${Z_{\chi}:=1/(q;q)_{\infty}}$ is the normalization constant appearing in \eqref{t11-3}, we get \eqref{t12-1}.
 Equivalence between \eqref{t12-1} and \eqref{t12-3} will be shown in Section~\ref{s23} although it is rather straightforward.
 The key observation is the formula
 \begin{align}
 \label{chiS}
 \P(\chi+S\le n)
 =
 \frac{1}{\bigl(-t q^{\frac12+n};q\bigr)_{\i}}.
 \end{align}
 This was already shown in \cite[Section 2]{BeteaBouttier2019PeriodicSchur}, but we will give a more direct proof in Lemma~\ref{l27} below. From \eqref{chiS}, we easily see that \eqref{t12-3} can be rewritten as{\samepage
 \[\P(\mu_1+\chi+S\le n)
 =
 \P(\l_1+S\le n),\]
 which is clearly equivalent to \eqref{t12-1}.}

 Now our goal is to prove \eqref{t12-3}. At this point we remark that, thanks to results from integrable probability \cite{Borodin2007PeriodicSchur,ImamuraSasamoto2019qTASEPtheta}, the expectation and the probability on both sides are known to be written as Fredholm determinants, see Propositions~\ref{p22} for the expectation on the left-hand side and Propositions~\ref{p25} for the probability on the right-hand side.
 Crucially, the kernels of the two Fredholm determinants turn out to be very similar: they both possess complex contour integral expressions with the integrands
 being the same and
 they differ from each other only by a choice of the contours.
 What we actually prove in this paper is that this difference of contours does not produce a difference to the value of the Fredholm determinants, see Theorem~\ref{t31}.\looseness=1

 \subsection{Background from integrable probability}

 In integrable probability, measures on the set of partitions written in terms of a pair of symmetric functions have attracted much attention. The simplest case, corresponding to the Cauchy identity \eqref{1-1}, is the Schur measure, introduced in~\cite{Okounkov2001InfiniteWedge}. For a partition $\l\in\Pe$, it is defined~by\looseness=1
 \begin{align}
 \label{1-15}
 M_{\text{S}}(\l)=s_{\l}(a_1,\dots,a_N)s_{\l}(b_1,\dots,b_M)/Z_{\text{S}},
 \end{align}
 where the normalization constant $Z_{\text{S}}$ is the right-hand side of \eqref{1-1}.
 From the fact that Schur polynomials can be written as single determinants due to the Jacobi--Trudi identity, the Schur measure defines a determinantal point process (DPP) \cite{Okounkov2001InfiniteWedge}, which means that all correlation functions are written as determinants with a common kernel, see, e.g., \cite{Soshnikov2000DPP}.

The Schur measure has connections with a few variants of the totally asymmetric simple exclusion process (TASEP) or the last passage percolation, which are typical models in the Kardar--Parisi--Zhang universality class. Such connection was first discovered by Johansson \cite{Johansson2000}, who showed that the current of the TASEP has the same law as the marginal $\l_1$ in the Schur measure.
 For our purposes a different variant of the TASEP with a certain pushing mechanism, called the \emph{PushTASEP}, is relevant
 \cite{Borodin2008,Warren-Windridge2009}.
 The precise definition of the dynamics can be found in the references.
 By using certain Markov dynamics on the Gelfand--Tsetlin cone,
 one can show that
 the PushTASEP is associated with the Schur measure. If $X_N(M)$ is the position of the~$N$-th particle at time $M$ starting from the so called \emph{step initial condition} $X_k(0)=k$ for~${k\in \mathbb{N}}$, we have the following relation:
 \begin{align}
 \label{1-16}
 \P (X_N(M)-N\le n)=\P(\l_1\le n),
 \end{align}
 where the right-hand side is the marginal of $\l_1$ from the Schur measure and can be written as
 a Fredholm determinant.

 The PushTASEP possesses a solvable $q$-deformation called $q$-PushTASEP, introduced in \mbox{\cite{BorodinPetrov2013NN, MatveevPetrov2015}}. This model is of particular interests as it interpolates between various notable models in the KPZ class including the log-gamma polymer model and the KPZ equation \cite{BoCoFeVe2015,MatveevPetrov2015}.
 Employing randomized versions of the RSK correspondence or of the Yang--Baxter equation one can relate the $q$-PushTASEP with the marginal $\mu_1$ of the $q$-Whittaker measure~\eqref{1-3} in the same way as in~\eqref{1-16} \cite{BufetovMucciconiPetrov2019,MatveevPetrov2015,MucciconiPetrov2020qWprocess}. Unlike in the $q=0$ case, such matching, although useful, does not lead immediately to explicit solutions of the model and in order to extract manageable formulas to study the $q$-PushTASEP additional considerations are needed. This is because, up to now, relations between the $q$-Whittaker measure and DPP had not sufficiently well understood. A~relevant result in this direction is a matching discovered by Borodin in~\cite{Borodin2018momentmatch} between certain multiplicative expectations of the Macdonald measure and Schur measure; we are going to comment on this in Section~\ref{sec:comparison}. Nevertheless, through the use of Macdonald difference operators, Markov duality, Bethe ansatz or elliptic determinants \cite{BorodinCorwin2014Mac,bcs2014,ImamuraSasamoto2019qTASEPtheta}, Fredholm determinant formulas describing interesting expectations of marginals of the $q$-Whittaker measure had been found. In general, following such approaches, one arrives at the following result:
 \begin{align}
 \label{qEdet}
 \E\left[\frac{1}{(\zeta q^{-\mu_1};q)_{\i}}\right]
 =
 \det(1- K)_{\ell^2(\Z)},
 \end{align}
 where $\mu_1$ is a marginal of the $q$-Whittaker measure, $\zeta$ is a parameter and $K$ is some kernel.

 In the references listed above several expressions for the right-hand side of \eqref{qEdet} have been obtained. Techniques involving Macdonald difference operators or Bethe ansatz lead to Fredholm determinant formulas where the kernel $K$ does not suggest any immediate relation with natural DPPs.
 On the other hand the elliptic determinant approach introduced in \cite{ImamuraSasamoto2019qTASEPtheta} produces a~substantially simpler formula where $K$ becomes a product between a finite rank kernel and the Fermi--Dirac distribution. Since such expressions, in particular the Fermi--Dirac distribution, are well known to appear in the description of free fermions at finite temperature, findings of \cite{ImamuraSasamoto2019qTASEPtheta} suggest a clear connections between discrete KPZ solvable models and DPP.

 This observation has led us to a careful comparison between the structure of the $q$-Whittaker measure and the (shift-mixed) periodic Schur measure \cite{Borodin2007PeriodicSchur}, as the latter is a canonical model of free fermions at positive temperature in one dimension \cite{BeteaBouttier2019PeriodicSchur}.
 In recent years the same model has received some attention and its properties and generalizations have been considered in \cite{Ahn-Russkikh-VanPeski2021GFFcylinder,BeteaBouttNejVul2018freeBoundarySchur,Koshida2021periodicmac}.
 The key observation we make is that the Fredholm determinant appearing in the shift-mixed periodic Schur measure is almost the same as the one found in \cite{ImamuraSasamoto2019qTASEPtheta} for the $q$-Whittaker measure, the only difference being in the contours of the contour integral expression for the kernels. In fact establishing this equivalence of the two Fredholm determinants is technically our main result as given in Theorem \ref{t31}.

 Before our result \eqref{t12-1}, outside of the TASEP/Schur case described above, explicit connections between KPZ solvable models and positive temperature free fermions
 had remained concealed. One notable exception to this claim was presented by the KPZ equation itself. The celebrated solution of the KPZ equation derived in 2010 by \cite{Amir-Corwin-Quastel2011KPZ,Calabrese-leDoussal-Rosso2010KPZ,Dotsenko2010KPZ,Sasamoto-Spohn2010KPZ-NPB,Sasamoto-Spohn2010KPZ-PRL} provides, in fact, a~clear relation between the one point function of the random KPZ height with narrow wedge initial condition and one point edge statistics of a system of trapped free fermions at finite temperature~\cite{Dean-LeDoussal-Majumdar-Schehr2015KPZfermion}. See also \cite{Johansson2007AiryFermi, Liechty-Wang2020fermion}.
 The identity \eqref{t12-1} extends this fact to a much higher level, widening substantially the bridge between solvable KPZ models and free fermions at finite temperature.

 An additional observation, relating this time the one point function of the KPZ equation to multiplicative statistics of the Airy point process was made in \cite{Borodin-Gorin2016momentmatchKPZ}. A more general instance of this finding was discovered by Borodin in \cite{Borodin2018momentmatch}, who proved a matching of the expectations between the certain stochastic vertex models and multiplicative statistics of the Schur measure. A specialization of such matching in distribution to the ASEP leads to a connection with a~variant of the Laguerre ensemble, explored in \cite{Borodin-Olshanski2017ASEP-DPP}.
 Results obtained in this stream of works and especially in \cite{Borodin2018momentmatch} are indeed related to the one given in Theorem~\ref{thm:intro_prob} and the connection will be explained in Section~\ref{sec:comparison}.

 All the matching discussed in the previous paragraphs, as the one obtained by us in this paper, consist essentially in comparisons between various formulas. In our upcoming paper \cite{ImamuraMucciconiSasamoto2021skewRSK}, we will further improve on this aspect developing a bijective theory providing a combinatorial explanation of identities between models.

\subsection{Outline}
 This paper is organized as follows. In Section~\ref{s2}, we give Fredholm determinant formulas for the $q$-Whittaker and shift-mixed Schur measures obtained respectively in \cite{ImamuraSasamoto2019qTASEPtheta} and \cite{Borodin2007PeriodicSchur}.
 We also discuss the equivalence between Theorems~\ref{t1} and~\ref{t12}. The proof of Theorem~\ref{t12} consists in the equivalence of these Fredholm determinant formulas and it is given in Section~\ref{s3}. A key passage to establish the proof is given by Lemma~\ref{l32}\,(i)--(iii). In Section~\ref{sec:comparison}, we compare Theorem~\ref{thm:intro_prob} with one of the main results in \cite{Borodin2018momentmatch}. In Appendix~\ref{aa}, we summarize some basic notions which we frequently use in this paper including partitions, $q$-series. In Appendix~\ref{ab}, we prove useful bounds which allow us to show Proposition~\ref{p33}\,(i) in Section~\ref{s33}. In Appendix~\ref{ac}, we state
 and show estimates of the $q$-Pochhammer symbol used in Section~\ref{s34}.

\section[Fredholm determinants for q-Whittaker and shift-mixed periodic Schur measures]{Fredholm determinants for $\boldsymbol{q}$-Whittaker\\ and shift-mixed periodic Schur measures}
\label{s2}

The main purpose of this section is to give Fredholm determinant formulas for both sides of~\eqref{t12-3}; see Propositions~\ref{p22} and~\ref{p25}.
In Section~\ref{s3}, we will show that the two Fredholm determinants are in fact the same,
as stated in Theorem~\ref{t31}. In Section~\ref{s23}, we also give a proof of equivalence between Theorems~\ref{t1} and~\ref{t12}.

Hereafter, we consider set of variables $a_1,\dots,a_N$, $b_1,\dots,b_N$, i.e., we set $M=N$ in the notation used in Section~\ref{s1}. The general case $M \neq N$ of each result is recovered by considerations similar to the one reported in Remark~\ref{rem:intro}. We will denote the maximum (resp.\ minimum) value of $a_i$ and $b_i$ $i=1,\dots,N$ by $a_{\max}$, $b_{\max}$ (resp.\ $a_{\min}$, $b_{\min}$). We also always omit the factor $1/\big(2\pi \sqrt{-1}\big)$ in contour integrals.

\subsection[Fredholm determinant formula for the q-Whittaker measure]{Fredholm determinant formula for the $\boldsymbol{q}$-Whittaker measure}
\label{s21}

We focus on the expectation value $\E[1/(\zeta q^{-\mu_1};q)_{\i}]$ with respect to the $q$-Whittaker measure defined by \eqref{s13}. This quantity is the $q$-Laplace transform of the marginal $\mu_1$ and it admits explicit Fredholm determinant expressions \cite{BorodinCorwin2014Mac,BoCoFeVe2015,BufetovMucciconiPetrov2019,ImamuraSasamoto2019qTASEPtheta,MatveevPetrov2015}. For our arguments, we will make use of the formula that first appeared in \cite{ImamuraSasamoto2019qTASEPtheta}, which is conceptually different from the ones of the other references.

The following proposition is obtained by applying a few results
in \cite{ImamuraSasamoto2019qTASEPtheta}.

\begin{figure}[ht]
 \centering
 \includegraphics[scale=1]{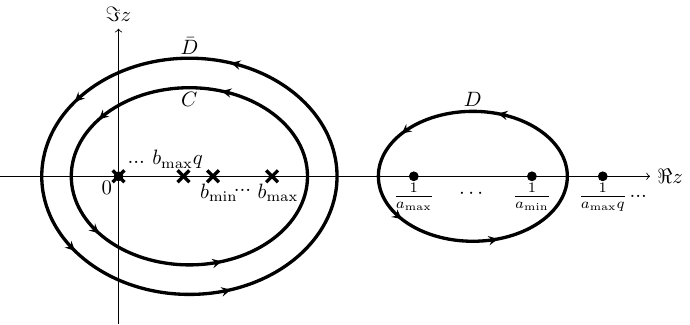}
 \caption{The contours $C$, $D$, and $\bar{D}$ appearing in \eqref{p222}
 and \eqref{t252}. ``$\bullet$''s and ``\boldcross''s represent the poles of
 $w$ and $z$ respectively.}
 \label{f1}
\end{figure}

\begin{Proposition}\label{p22}
Let $\zeta\in\C\setminus\{q^n\}_{n\in\Z}$, and
$a_i, b_i\in\R_{>0}$, $i=1,\dots,N$ satisfy $a_{\max}b_{\max}<1$ and~${qa_{\max}<a_{\min}}$.
Let $\mu_1$ be the first row of $\mu$ distributed according to the $q$-Whittaker measure~\eqref{s13}.
We have
\begin{align}
\label{p221}
\E\left[\frac{1}{(\zeta q^{-\mu_1};q)_{\i}}\right]
=
\det(1-f_{\zeta} K)_{\ell^2(\Z)},
\end{align}
where
$f_{\zeta}(m)=\frac{-\zeta q^m}{1-\zeta q^m}$
and
\begin{align}
\label{p222}
&K(m_1,m_2)=
\int_C\frac{{\rm d}z}{z}\int_{D}\frac{{\rm d}w}{w}
g_{a,b}(z,w;m_1,m_2),
\\
&
g_{a,b}(z,w;m_1,m_2)
=
\frac{w^{m_2}}{z^{m_1}}
\frac{w}{z-w}
\prod_{i=1}^N
\frac{(a_i z;q)_{\i}}{(a_iw;q)_{\i}}
\frac{(b_i/w;q)_{\i}}{(b_i/z;q)_{\i}},
\label{p223}
\end{align}
where the contour $C$ encloses anticlockwise
the poles at $z=0, b_iq^j$, $i=1,\dots,N$,
$j\in\Z_{\ge 0}$
whereas~$D$ encloses anticlockwise
the poles at
$w=1/a_i$, $i=1,\dots,N$
but not the ones at $w=0, 1/a_iq^j$,
$i=1,\dots,N$, $j\in\Z_{\ge 1}$.
They are illustrated in Figure~{\rm\ref{f1}}.
\end{Proposition}

To prove Proposition~\ref{p22}, we recall the discussions in \cite[Section 3]{ImamuraSasamoto2019qTASEPtheta}. There we consider a~two-sided $q$-Whittaker measure, which is a~probability measure on the set of signatures,
\begin{align*}
 \S_N=\big\{\l=(\l_1,\dots,\l_N)\in \Z^N\mid \l_1\ge\cdots\ge \l_N\big\}.
\end{align*}
Note that $\l_i$, $i=1,\dots,N$ can be negative. Such probability measures on the extended space can be realized by replacing the function $Q_{\mu}(x)$ in \eqref{s13} with $Q_{\mu}(\rho)$, where $\rho$ is a certain specialization of the algebra of symmetric functions.
Here we state only a particular case which corresponds to \cite[Definition 3.4]{ImamuraSasamoto2019qTASEPtheta}, setting $a_i=x_i$, $i=1,\dots,N$, $t=0$ and $\a_j=b_j$, $j=1,\dots,N$. With this choice $M_{q\text{W}}(\mu)$ in \eqref{s13} becomes
\begin{align}
\label{p224}
P_{\l}(x_1,\dots,x_N)Q_{\l}(\rho^{-}_b)/\Pi^{(-)}(x;b),
\end{align}
where $\l\in\S_N$, the parameters $x_i$, $b_i$,
satisfy $b_{\max}/x_{\min}<1$ and
\begin{align}
Q_\l(\rho^-_b)
 =\prod_{i=1}^{N-1}\big(q^{\l_i-\l_{i+1}+1};q\big)_{\infty}
\int_{C^N}\prod_{i=1}^N\frac{{\rm d}z_i}{z_i}\cdot
P_\l (1/z)
 \Pi^{(-)}(z;b)m_N^q(z).
\label{defQ}
\end{align}
Here $C$ is the same contour as the
one in \eqref{p222}, i.e., it encloses
positively
the poles at $z=0, b_iq^j$, $i=1,\dots,N$
and $j\in\Z_{\ge 0}$, $1/z$ in $P_\l$ is a shorthand notation for
$(1/z_1,\dots,1/z_N)$ and
\begin{align*}
&m_N^q(z)=\frac{1}{N!}\prod_{1\le i<j\le N}(z_i/z_j;q)_{\infty}(z_j/z_i;q)_{\infty},
\qquad
\Pi^{(-)}(z;b)=\prod_{i,j=1}^N\frac{1}{(b_i/z_j;q)_{\infty}}.
\end{align*}
For more detail see \cite[Appendix B]{ImamuraSasamoto2019qTASEPtheta}.

From \cite[Lemma 3.3, equation~(3.15)]{ImamuraSasamoto2019qTASEPtheta}, one sees that $Q_{\l}(\rho_b^-)=0$ unless
$0\ge \l_1\ge\cdots \ge \l_N$ and thus the measure
\eqref{p224} has support only on the negative region. Although the measures~\eqref{s13} and \eqref{p224} have different support, they enjoy the following relation:
\begin{Lemma} \label{l24}
Set $a_i=1/x_i$, $i=1,\dots,N$ and define the signature $\mu$ by the negation of $\l$,
$\mu:=(-\l_N,-\l_{N-1},\dots,-\l_1)$.
 Then we have
\begin{align}
\label{2-12d}
 P_{\l}(x_1,\dots,x_N)Q_{\l}(\rho^{-}_b)/\Pi^{(-)}(x;b)
=
M_{q{\rm W}}(\mu),
 \end{align}
 and in particular
\begin{align}
\label{2-12}
 \E_{q{\rm W}(x;\rho_b^-)}\left[\frac{1}{\big(\zeta q^{\l_N};q\big)_{\i}}\right]
 =
 \E_{q{\rm W}(a;b)}\left[\frac{1}{(\zeta q^{-\mu_1};q)_{\i}}\right],
\end{align}
 where \smash{$\E_{q{\rm W}(x;\rho_b^-)}[\cdot]$}
 and $\E_{q{\rm W}(a;b)}[\cdot]$
 are averages over the measures
 \eqref{p224} and \eqref{s13}, respectively.
\end{Lemma}
\begin{proof}
 Equation \eqref{2-12d} follows straightforwardly from a combination of definition \eqref{p224}
 and relations \cite[equations~(3.14) and (3.15)]{ImamuraSasamoto2019qTASEPtheta},
 \[
 P_{\l}(x_1,\dots,x_N)=P_{\mu} (a_1,\dots,a_N),
 \qquad
 Q_{\l}(\rho^-_b)=Q_{\mu}(b_1,\dots,b_N),
 \]
 once we recall hypothesis relating $a_i$ (resp.\ $\mu$)
 to $x_i$ (resp.\ $\l$).
\end{proof}

\begin{proof}[Proof of Proposition~\ref{p22}]
 It follows immediately from Lemma~\ref{l24}
 and the Fredholm determinant formula
 for the left-hand side of \eqref{2-12}.
 Changing the notations $a_i\rightarrow x_i$,
 $\alpha_j\rightarrow b_j$,
 $v\rightarrow w$,
 $n\rightarrow -m_1$, $m\rightarrow -m_2$
 and setting $N=M$ and $t=0$
 in~\cite[equation~(4.37)]{ImamuraSasamoto2019qTASEPtheta},
 we obtain the Fredholm determinant formula
\begin{align}
\label{p225}
\E_{q\text{W}(x;\rho_b^-)}\left[\frac{1}{\big(\zeta q^{\l_N};q\big)_{\i}}\right]
=
\det(1-f_{\zeta} M)_{\ell^2(\Z)}.
\end{align}
Here $f_{\zeta}(m)$ is defined above~\eqref{p222}
and
\begin{align*}
M(m_1,m_2)=
\int_C\frac{{\rm d}z}{z}\int_{D}\frac{{\rm d}w}{w}
\frac{w^{m_2}}{z^{m_1}}
\frac{w}{z-w}
\prod_{i=1}^N
\frac{(z/x_i;q)_{\i}}{(w/x_i;q)_{\i}}
\frac{(b_i/w;q)_{\i}}{(b_i/z;q)_{\i}},
\end{align*}
where $C$ is the same contour as \eqref{defQ} and $D$ encloses positively
$x_i,~i=1,\dots,N$.
Combining~\eqref{2-12} with~\eqref{p225},
we get~\eqref{p221}.
\end{proof}

\subsection[Fredholm determinant formula for the shift-mixed periodic Schur measure]{Fredholm determinant formula \\
for the shift-mixed periodic Schur measure}\label{s22}

We have discussed in Section~\ref{s1} how the periodic Schur measure $M_{\text{pS}}(\l)$ \eqref{1-8} becomes a determinantal point process once all the rows of $\lambda$ get shifted by an independent random quantity $S$ distributed as \eqref{1-11}. This was originally observed by Borodin in \cite{Borodin2007PeriodicSchur} and the measure describing the process
\begin{align}
\label{216}
(\l_1+S,\l_2+S,\dots)
\end{align}
takes the name of \emph{shift-mixed periodic Schur measure}. The probability distribution of the first element $\l_1+S$ can be represented as a single Fredholm determinant and this is reported in the next proposition.

\begin{Proposition}\label{p25}
 Let $\l_1+S$ be the first element of \eqref{216} distributed according to the shift-mixed periodic Schur measure.
 For $k\in\Z_{\ge 0}$, we have
\begin{align} \label{p251}
 \P (\l_1+S\le k )
 =\det (1+f L )_{\ell^2(\Z)},
\end{align}
 where
 $f(m)$ is given in terms of
 $f_{\zeta}(m)$ defined above \eqref{p222},
\begin{align}
\label{p251-1}
 f(m)=f_{\zeta}(m)|_{\zeta=-tq^{1/2+k}}
 =\frac{tq^{1/2+k+m}}{1+tq^{1/2+k+m}}
\end{align}
and
by using $g_{a,b}(z,w,m_1,m_2)$~\eqref{p223},
the kernel
$L(m_1,m_2)$ is defined by
\begin{align}
&L(m_1,m_2)=
\int_C\frac{{\rm d}z}{z}\int_{\bar{D}}\frac{{\rm d}w}{w}
g_{a,b}(z,w,m_1,m_2).
\label{t252}
\end{align}
Here contour $C$ is oriented counterclockwise and it encloses the origin with its radius
$r$ satisfying~${b_{\max}<r<1/a_{\max}}$
whereas $\bar{D}$ encloses $C$
anticlockwise.
$C$ and $\bar{D}$ are depicted in Figure~{\rm\ref{f1}}
given below Proposition~{\rm\ref{p22}}.

\end{Proposition}

\begin{proof}
From \cite[Theorem 2.2]{Borodin2007PeriodicSchur}, we have, for $k\in\Z_{\ge 0}$,
\begin{align*}
\P(\l_1+S\le k)
=\det\big(1-\tilde{L}\big)_{\ell^2(k+1,k+2,\dots)},
\end{align*}
where the kernel $\tilde{L}$ is defined as
\begin{align}
\tilde{L}(\ell_1,\ell_2)={}&\int_{|w|=r}\frac{{\rm d}w}{w}\nonumber\\
&\times\int_{|z|=r'}\frac{{\rm d}z}{z}\frac{z^{\ell_2-1}}{w^{\ell_1-1}}
\prod_{i=1}^N
\frac{(a_iz;q)_\i}{(a_i w;q)_{\i}}
\frac{(b_i/w;q)_\i}{(b_i/z;q)_\i}
\cdot
\sum_{m\in \Z}
\frac{t q^{\frac12+m}}{1+t q^{\frac12+m}}
\left(\frac{w}{z}\right)^m,
\label{112}
\end{align}
and $r$, $r'$ satisfy\footnote{In~\cite{BeteaBouttier2019PeriodicSchur,Borodin2007PeriodicSchur},
the conditions for $r$, $r'$ are
given by using the analyticity
of the integrand for the periodic
Schur process with positive specializations.
This condition can be translated into
the explicit one \eqref{t111}
in our case~\eqref{1-8}.
}
\begin{align}
b_{\max}<r,r'<\frac{1}{a_{\max}},\qquad 1<\frac{r}{r'}<q^{-1}.
\label{t111}
\end{align}
We will now show that $\det\big(1-\tilde{L}\big)_{\ell^2(k+1,k+2,\dots)} =\det(1+fL)_{\ell^2(\Z)}$. Note that the kernel $\tilde{L}(\ell_1,\ell_2)$ \eqref{112} can be written as the kernels of the operators $A$ and $B$, which map between $(k+1,k+2,\dots)$ and $\mathbb{Z}$:
\begin{align*}
\tilde{L}(\ell_1,\ell_2)=\sum_{m\in\Z}A(\ell_1,m)
B(m, \ell_2),
\end{align*}
where
\begin{align*}
A(\ell,m)&=\int_{|w|=r}\frac{{\rm d}w}{w} w^{m-\ell+1}
\prod_{i=1}^N
\frac{(b_i/w;q)_\i}{(a_i w;q)_{\i}},
\\
B(m,\ell)&=
\frac{t q^{1/2+m}}{1+t q^{1/2+m}}
\int_{|z|=r'}\frac{{\rm d}z}{z} z^{\ell-1-m}
\prod_{i=1}^N
\frac{(a_iz;q)_\i}{(b_i/z;q)_\i}.
\end{align*}
By a basic property of the determinant, we have
\begin{align*}
\det(1-AB)_{\ell^2(k+1,k+2,\dots)}=\det(1-BA)_{\ell^2(\Z)},
\end{align*}
and evaluating the kernel in the right-hand side, we find
\begin{gather*}
(BA)(n_1,n_2)
\\
\qquad=
\frac{t q^{n_1+\frac12}}{1+t q^{n_1+\frac12}}
\int_{|z|=r'}\frac{{\rm d}z}{z}
\int_{|w|=r}\frac{{\rm d}w}{w}
\frac{w^{n_2}}{z^{n_1}}
\prod_{i=1}^N
\frac{(a_iz;q)_\i}{(a_i w;q)_{\i}}
\frac{(b_i/w;q)_\i}{(b_i/z;q)_{\i}}
\cdot
\sum_{\ell=k+1}^{\i}\left(\frac{z}{w}\right)^{\ell-1}
\\
\qquad=
\frac{-t q^{\frac12+k+m_1}}{1+t q^{\frac12+k+m_1}}
\int_{|z|=r'}\frac{{\rm d}z}{z}
\int_{|w|=r}\frac{{\rm d}w}{w}
\frac{w^{m_2}}{z^{m_1}}
\prod_{i=1}^N
\frac{(a_iz;q)_\i}{(a_i w;q)_{\i}}
\frac{(b_i/w;q)_\i}{(b_i/z;q)_{\i}}
\cdot
\frac{w}{z-w}
\\
\qquad=
-f(m_1)L(m_1,m_2).
\end{gather*}
Here in the second equality we set $m_i=n_i-k$ for $i=1,2$.
\end{proof}

\subsection[Equivalence of Theorems 1.1 and 1.3]{Equivalence of Theorems \ref{t1} and \ref{t12}}
\label{s23}

This subsection is devoted to proving that Theorem~\ref{t1} is equivalent to Theorem~\ref{t12}. Since equivalence between Theorems~\ref{t1} and~\ref{t12}\,(a) is straightforward, the only nontrivial part is to show the equivalence between Theorem~\ref{t12}\,(a) and (b), i.e., \eqref{t12-1} and \eqref{t12-3}. This is a~consequence of the computation reported in the following lemma.
\begin{Lemma}
\label{l27}
For $n\in\Z$, $0<q<1$ and
$t>0$, we have
 \begin{align}
\label{2-31}
 \P(\chi+S\le n)
 =
 \frac{1}{\bigl(-t q^{\frac12+n};q\bigr)_{\i}},
\end{align}
 where $\chi$ and $S$ are independent random variables with distributions
 \eqref{t11-3} and \eqref{1-11}, respectively.
\end{Lemma}
 \begin{proof}
 We write the left-hand side of \eqref{2-31} as
 \begin{align}
 \label{2-34}
 \P(\chi+S\le n)
 =
 \sum_{\ell\in\Z} \P(S=\ell)\P(\chi\le n-\ell).
 \end{align}
 Applying \eqref{1-11} and
 \begin{align*}
 \P(\chi\le m)=
 \begin{cases}
 (q;q)_{\i}/(q;q)_m,& m\in\Z_{\ge 0},
 \\
 0,& m\in\Z_{< 0},
 \end{cases}
 \end{align*}
 which follows from \eqref{t11-3} and the fact $\sum_{n=0}^mq^n/(q;q)_n=1/(q;q)_m$, $m=0,1,2,\dots$, which can be checked easily by mathematical induction, to the right-hand side of \eqref{2-34}, we find
 \begin{align*}
 \P(\chi+S\le n)
 =
 \frac{1}{\theta\bigl(-q^{1/2}t\bigr)}\sum_{\ell=-\infty}^n
 \frac{t^{\ell}q^{\frac{\ell^2}{2}}}
 {(q;q)_{n-\ell}}
 =
 \frac{t^nq^{\frac{n^2}{2}}\bigl(-q^{\frac12-n}t^{-1};q\bigr)_{\i}}{\theta\bigl(-q^{1/2}t\bigr)},
 \end{align*}
 where, in the second equality we used a version of the $q$-binomial theorem \eqref{a2-4}.
 Noting from the definitions of the $q$-Pochhammer symbol \eqref{a2-1} and \eqref{a2-1d}, the numerator of the last expression becomes
 \begin{align*}
 t^n q^{n^2/2} \bigl(-q^{\frac12-n}t^{-1};q\bigr)_\i
=\bigl(-t q^{1/2};q\bigr)_n \bigl(-q^{1/2}/t;q\bigr)_{\i},
 \end{align*}
 we arrive at the right-hand side of \eqref{2-31}.
 \end{proof}

 By using \eqref{2-31}, we see that
 the left-hand side of \eqref{t12-3}
 is expressed as
 \begin{align*}
 \E
 \left[ \frac{1}{\bigl(-t q^{\frac12+n-\mu_1}\bigr)}\right]
 =\P(\mu_1+\chi+S\le n).
 \end{align*}
 Thus \eqref{t12-3} can be rewritten as
 \begin{align*}
 \P(\mu_1+\chi+S\le n)
 =
 \P(\l_1+S\le n),
 \end{align*}
 which is clearly equivalent to \eqref{t12-1}.

\section[Equivalence of two Fredholm determinants by shift of contours]{Equivalence of two Fredholm determinants\\ by shift of contours} \label{s3}

In this section, we will force some conditions on parameters $a_i$, $b_i$,
$i=1,\dots,N$. Namely we assume that $a_i,b_i\in (0;1)$, $i=1,\dots,N $ are distinct and ordered as $a_1>a_2>\cdots > a_N$ (hence~${a_{\max}=a_1}$ and $a_{\min}=a_N$)
and moreover
\[
 a_{1}/a_{N}<q^{-\frac12+\e},
\]
for a fixed $\e \in (0,1/2)$. Recall that for the sake of our main result Theorem~\ref{t1} these restrictions do not constitute a problem, see Remark~\ref{rem:intro}.

\subsection{Restating the main theorem as a determinant identity}
\label{s31}
In the previous section, precisely in Propositions~\ref{p22} and~\ref{p25}, we have seen that both sides of~\eqref{t12-3} can be written as Fredholm determinants. By this equality \eqref{t12-3} is equivalent to the following identity for the two Fredholm determinants.
\begin{Theorem}
\label{t31}
Let $K(m_1,m_2)$ be the kernel given in ~\eqref{p222} and $L(m_1,m_2)$ in ~\eqref{t252}.
Also let~$f$ be given by \eqref{p251-1}.
Then we have
\begin{align}
\label{31}
\det(1-f K)_{\ell^2(\Z)}
=
\det(1+fL)_{\ell^2(\Z)}.
\end{align}
\end{Theorem}
In this section, we will elaborate the proof of Theorem~\ref{t31}.
Comparing the expressions of the two kernels in~\eqref{p222} and in~\eqref{t252}, we notice that both of them are written in the double integral forms with the common integrand $g_{a,b}(z,w;m_1,m_2)$, given by \eqref{p223}. The only difference between them is the integration contour of $w$: the contour $D$ in $K(m_1,m_2)$ encloses $1/a_j$, $j=1,\dots,N$, while $\bar{D}$ in $L(m_1,m_2)$ encloses the origin and the contour of $z$. See Figure~\ref{f1}.

Result of Theorem~\ref{t31} is not trivial outside of the case $q=0$, when the kernel themselves are equal, i.e., $K(m_1,m_2)=-L(m_1,m_2)$ for all $m_1,m_2\in\Z_{\le 0}$. This can be shown as follows. When $q=0$, we notice that $f(m)$ becomes the indicator function of $\mathbb{Z}_{\le 0}$ and the integrand~${g_{a,b}(z,w;m_1,m_2)}$ in \eqref{p223} becomes
\begin{align*}
 g_{a,b}(z,w;m_1,m_2)
 =\frac{w^{m_2}}{z^{m_1}}
\frac{w}{z-w}
\prod_{i=1}^N
\frac{(1-a_i z)}{(1-a_iw)}
\frac{(1-b_i/w)}{(1-b_i/z)},
\end{align*}
and we see that the poles in the $w$ variable include only $w=0, z, 1/a_i$, $i=1,\dots,N$ and in particular $\infty$ is not a pole since $m_2 \le 0$. Thus using Cauchy's integral theorem, we see that the contour $D$ can be deformed to $\bar{D}$ without changing the value of the kernel up to sign.

For general $0<q<1$, this is no longer the case. In this case, equality does not hold for the kernels, but it only does at the level of the Fredholm determinants \eqref{31}. To prove this we consider a sequence of changes of the contours of $w$ starting from $D$, and including, sequentially, poles of the form $1/a_j q^k$, $j=1,\dots,N$ starting from $k=0,1,\dots$ and finally including the singularity at infinity. We find that such modifications of the kernel do not affect its Fredholm determinant, so that, once all the singularities have been added in the contour one can switch such contour to $\bar{D}$ by using the Cauchy theorem.

\begin{figure}[ht]
 \centering
 \includegraphics[scale=1]{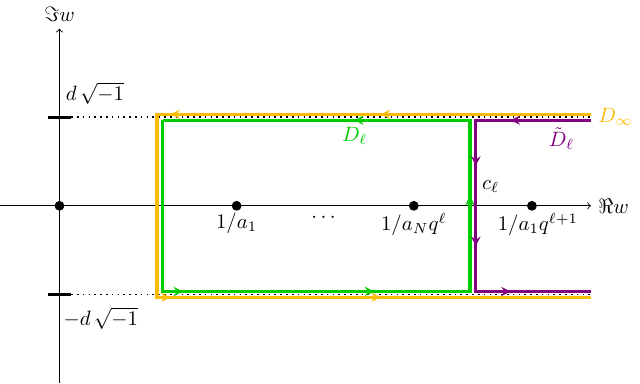}
 \caption{Contours $D_{\ell}$, $\tilde{D}_\ell$, and $D_{\i}$
 are illustrated.
 We set $c_{\ell}$
 to be the center of the interval $\big[\frac{1}{a_{N}q^\ell},
 \frac{1}{a_{1} q^{\ell+1}}\big]$, i.e., $c_{\ell}=\frac{1}{2q^\ell}\big(\frac{1}{a_{1}q}+\frac{1}{a_{N}}
 \big)$.}
 \label{fig2}
\end{figure}

Notations for the series of contours just described, their corresponding kernels and Fredholm determinants are provided next.
\begin{Definition}\label{d31}
For $\ell\in\mathbb{Z}_{\ge 0}$, we define the kernels
\begin{align}
\label{d31-3}
&K_\ell (m_1,m_2)=
\int_C\frac{{\rm d}z}{z}\int_{D_\ell}\frac{{\rm d}w}{w} g_{a,b}(z,w;m_1,m_2),
\\
\label{d31-4}
&K_{\i}(m_1,m_2)
=
\int_C\frac{{\rm d}z}{z}\int_{D_{\i}}\frac{{\rm d}w}{w} g_{a,b}(z,w;m_1,m_2).
\end{align}
The function $g_{a,b}$ is given by \eqref{p223} and the contour $C$ encloses the origin and it has radius larger than $b_{\max}$ as in \eqref{p222}, \eqref{t252}. The integration contours $D_\ell$ and $D_{\infty}$ are illustrated in Figure~\ref{fig2}.
Note $D_\ell$ encircles poles $1/q^{i}a_j$, $i=0,\dots,\ell$, $j=1,\dots,N$ and no other singularities lie inside~$D_{\ell}$.
Recalling the function $f(m)$ defined by \eqref{p251-1}, we further define the Fredholm determinants
\begin{align}
\label{d31-1}
&F_\ell:=\det(1-fK_\ell)_{\ell^2(\Z)}\qquad \text{for}\quad \ell=0,1,2,\dots,
\\
\label{d31-2}
&F_{\infty}:=\det(1-f K_{\i})_{\ell^2(\Z)}.
\end{align}
\end{Definition}

Note that by definition, we see $K_0(m_1,m_2)=K(m_1,m_2)$. Thus by Proposition~\ref{p222}, we immediately see $F_0=\det(1-fK)_{\ell^2(\Z)}$. In \cite[Proposition 4.7]{ImamuraSasamoto2019qTASEPtheta}, we showed that the kernel~${K(m_1,m_2)}$ is trace-class, thus $\det(1-fK)_{\ell^2(\Z)}$ is well defined.
This fact will be generalized to $F_{\ell}$ for any~${\ell\in\Z_{\ge 0}}$ and also $F_{\i}$ in Lemma~\ref{l32}\,(iii).

For this purpose and also for later use,
we will give different representations for
the kernels~\eqref{d31-3} and \eqref{d31-4}.
 To simplify notation, we introduce a relabelling of all the inverse of poles~$a_i q^j$ in
decreasing order and call them $\tilde{a}_r,r\in\Z_{>0}$. Namely, we write
\begin{align}
\label{s32-1}
 \tilde{a}_{r}=a_{k}q^u,
\end{align}
where $r=uN+k$ with $r\in\Z_{>0}$, $k\in\{1,\dots,N\}$ and $u\in\Z_{\ge 0}$.
Since we assumed $a_1>a_2>\cdots >a_N$ without loss of generality, $\tilde{a}_1>\tilde{a}_2>\cdots$ holds.

For $m\in\Z$ and $r\in\Z_{\ge 0}$,
we define two infinite-dimensional matrices,
$A(m;r)$, $B(r;m)$, as
\begin{align}
 \label{3-15d}
 A(m;r)&=f(m)\int_C\frac{{\rm d}z}{z^{1+m}}
 \frac{1}{z-\tilde{a}_{r}^{-1}}
 \prod_{i=1}^N\frac{(a_i z;q)_{\i}}{(b_i/z;q)_{\i}},
 \\
 \label{3-16d}
 B(r;m)&=\tilde{a}_{r}^{-m}
 \underset{w=\tilde{a}_{r}^{-1}}{\text{Res}}
 \prod_{i=1}^N\frac{(b_i/w;q)_{\i}}{(a_i w;q)_{\i}}.
\end{align}
For $\ell\in\Z_{\ge 0}$, we also introduce their truncated version, $A_{\ell}(m;r)$, $B_{\ell}(r;m)$, defined by
\[
 A_{\ell}(m;r):=A(m;r) \mathbf{1}_{\{1,\dots, N(\ell+1)\}}(r),\qquad
 B_{\ell}(r;m):= \mathbf{1}_{\{1,\dots, N(\ell+1)\}}(r)B(r;m),
\]
where for a set $A\subset\Z_{\ge0}$, $\mathbf{1}_A$ is the indicator function of $A$. By evaluating the poles inside $D_\ell$ of~\eqref{d31-3}, we easily find
\begin{align}
\label{3-10}
 &f(m_1)K_{\ell}(m_1,m_2)=\sum_{r=1}^{\infty}A_{\ell}(m_1;r)B_{\ell}(r;m_2).
\end{align}
Below in the proof of Proposition~\ref{p33}\,(iii), we will also see that
\begin{align}
\label{3-11}
 f(m_1)K_{\infty}(m_1,m_2)=\sum_{r=1}^{\i}A(m_1;r)B(r;m_2).
\end{align}
In the following discussions, it is also useful to consider the conjugated matrices
\begin{align}
 \label{3-16dd}
 \tilde{A}(m;r)&=\t(m) A(m;r)\sigma(r),
 \qquad
 \tilde{B}(r;m)=\sigma^{-1}(r) B(r;m)\t^{-1}(m),
\end{align}
where
\begin{align}
\label{3-12d}
 \t(m)=
 \begin{cases}
 a_{1}^{-m} q^{-\frac{m^2}{2N}-\frac{m}{2}+\e m},
 & m\ge 0,\\
 1,& m<0,
 \end{cases}
 \qquad
 \sigma(r)= q^{-(1-\omega)u}
\end{align}
for $\e\in(0,1/2)$ and $\o\in(0,1/2-\e)$, and also
their truncated versions $\tilde{A}_{\ell}(m;r)$, and
$\tilde{B}_{\ell}(r;m)$,
\begin{align}
\label{3-15at}
 \tilde{A}_{\ell}(m;r):=\tilde{A}(m;r) \mathbf{1}_{\{1,\dots, N(\ell+1)\}}(r),\qquad
 \tilde{B}_{\ell}(r;m):= \mathbf{1}_{\{1,\dots, N(\ell+1)\}}(r)\tilde{B}(r;m).
\end{align}
Note that the value of the
Fredholm determinant would not change when we replace
$A_\ell$, $B_\ell$,~$A$,~$B$ by their conjugated versions,
$\tilde{A}_\ell$, $\tilde{B}_\ell$, $\tilde{A}$, $\tilde{B}$
in the kernel \eqref{3-10}, \eqref{3-11}.

\begin{Proposition}\label{p33}
Assume $a_i,b_i\in(0,1)$, $i=1,\dots,N$, $a_1>\cdots>a_N$, and $a_{1}/a_{N}<q^{-\frac12+\e}$ with $\e\in(0,1/2)$.
\begin{itemize}\itemsep=0pt
 \item[{\rm (i)}]
The operators with
the kernels, $\tilde{A}(m;r)$, $\tilde{B}(r;m)$,
$\tilde{A}_\ell(m;r)$, $\tilde{B}_\ell(r;m)$
for $m\in\Z$ and $r,\ell\in\Z_{\ge 0}$
are of Hilbert--Schmidt class.
 \item[{\rm (ii)}] $K_{\i}(m_1,m_2)$
 can be expressed as the limit $\ell\rightarrow\i$ of $K_{\ell}(m_1,m_2)$,
 \begin{equation}
 \label{p33-1}
 K_{\i}(m_1,m_2)=\lim_{\ell\rightarrow\i}K_\ell(m_1,m_2).
 \end{equation}
 \item[{\rm (iii)}]
 The Fredholm determinants $F_\ell $ and $F_{\infty}$, defined by~\eqref{d31-1} and~\eqref{d31-2}, respectively, are well defined.
\end{itemize}
\end{Proposition}
\begin{proof}
 (i) By the estimates of $\tilde{A}(m;r)$, $\tilde{B}(r;m)$ in Corollary~\ref{cb2},
 it is easy to check
 the conditions of the Hilbert--Schmidt operators
\[
 \sum_{m\in\Z}\sum_{r=1}^{\i}\big|\tilde{A}(m;r)\big|^2<\i,
 \qquad
 \sum_{m\in\Z}\sum_{r=1}^{\i}\big|\tilde{B}(r;m)\big|^2<\i.
\]
Obviously, the same relations also hold for $\tilde{A}_{\ell}$, $\tilde{B}_{\ell}$.

 For (ii), let us introduce the notation,
\begin{align}
\label{3-26}
 K(m_1,m_2;\sharp)=
 \int_C\frac{{\rm d}z}{z}
 \int_{\sharp}\frac{{\rm d}w}{w}
 g_{a,b}(z,w,m_1,m_2), \qquad m_1,m_2\in\Z,
\end{align}
 where $g_{a,b}(z,w,m_1,m_2)$ is defined by \eqref{p223} and $\sharp\in\big\{D_{\ell}, \tilde{D}_{\ell}\big\}_{\ell=0,1,\dots}\cup\{D_{\i}\}$ represents a contour for $w$. They are illustrated in Figure~\ref{fig2}. By definition, one sees
\[
 K_\ell(m_1,m_2)=K(m_1,m_2;D_\ell),\qquad
 K_{\i}(m_1,m_2)= K(m_1,m_2;D_\i),
\]
Using~\eqref{3-26}, we show \eqref{p33-1}.
 Noting $D_{\i}=D_\ell\cup\tilde{D}_\ell$, we see that
\begin{align*}
 K(m_1,m_2;D_{\i})
 &=
 K(m_1,m_2;D_\ell)+ K\big(m_1,m_2;\tilde{D}_\ell\big)
 \\
 &=
 \lim_{\ell\rightarrow\i}K(m_1,m_2;D_\ell)
 +
 \lim_{\ell\rightarrow\i}K\big(m_1,m_2;\tilde{D}_\ell\big)
 =\lim_{\ell\rightarrow\i}K(m_1,m_2;D_\ell).
\end{align*}
 Here in the second equality, we used the
 fact that $D_\ell\cup\tilde{D}_\ell$ $(=D_{\i})$
 does not depend on $\ell$ while in the third equality we used the
 fact \smash{$\lim_{\ell\rightarrow\i}K\big(m_1,m_2;\tilde{D}_\ell\big)=0$}. This follows from estimates of Lemma \ref{lc1}\,(iii): the contour $\tilde{D}_\ell$ can be divided into three
 parts \smash{$\tilde{D}_\ell=\tilde{D}_\ell^{(1)}\cup \tilde{D}_\ell^{(2)}\cup \tilde{D}_\ell^{(3)}$}, where~\smash{$\tilde{D}_\ell^{(1)}$} \big(resp.\ \smash{$\tilde{D}_\ell^{(3)}$}\big) are horizontal half line from~${+\infty+ d\sqrt{-1}}$ to $c_{\ell}+d\sqrt{-1}$ (resp.\ from $c_{\ell}-d\sqrt{-1}$ to $+\infty-d\sqrt{-1}$) while \smash{$\tilde{D}_\ell^{(2)}$} is vertical line segment from $c_{\ell}+d\sqrt{-1}$ to $c_{\ell}-d\sqrt{-1}$. We get $\lim_{\ell\rightarrow\i}K\big(m_1,m_2;\tilde{D}_\ell\big)=0$ applying \eqref{lc1-3} for the cases \smash{$\tilde{D}_{\ell}^{(1)}$} and \smash{$\tilde{D}_{\ell}^{(3)}$} and \eqref{lc1-3'} for the case \smash{$\tilde{D}_{\ell}^{(2)}$} to the part \smash{$\prod_{i=1}^N(a_iw;q)_{\infty}^{-1}$} in the integrand \eqref{p223}.

(iii)
 By equations~\eqref{3-10}, \eqref{3-15at} and (i), we see that the operators with the kernel
 \[\t(m_1)f(m_1) K_\ell(m_1,m_2)\t^{-1}(m_2),\]
 $\ell=0,1,2,\dots$ are trace-class. Since the factor $\t(m_1)/\t(m_2)$ does not affect the value of the determinant, then $F_\ell$ in \eqref{d31-1} is well defined.
 Next combining \eqref{3-10} and~\eqref{p33-1}, we have \eqref{3-11}. Then by the same arguments used above, the operator with the kernel $\t(m_1)K_\i (m_1,m_2)\t^{-1}(m_2)$ is trace-class and $F_\infty$ in \eqref{d31-2} is well defined.
\end{proof}

Theorem~\ref{t31} can be immediately proved by the following lemma.

\begin{Lemma}\label{l32}
Denote the left-hand side of ~\eqref{31} by $F$, i.e., $F:=\det\left(1-f K\right)_{\ell^2(\Z)}$. The following statements hold.
\begin{itemize}\itemsep=0pt
\item[{\rm (i)}]
$F_\ell$ does not depend on $\ell=0,1,2,\dots$ and is equal to $F$, i.e., we have $F_{\ell}=F$, for any~${\ell\in\Z_{\ge 0}}$.
\item[{\rm (ii)}] We have $\lim_{\ell \rightarrow\i} F_{\ell}=F_{\i}$.
\item[{\rm (iii)}] We have $K_{\i}(m_1,m_2)=-L(m_1,m_2)$
and thus $F_{\i}=\det(1+f L)_{\ell^2(\Z)}$.
\end{itemize}
\end{Lemma}
The proofs of the above (i)--(iii) will
be given in Sections \ref{s32}--\ref{s34}, respectively.

\begin{proof}[Proof of Theorem~\ref{t31}]
 This is a straightforward consequence of Lemma~\ref{l32}, which implies
 \begin{align*}
 \det(1-fK)_{\ell^2(\Z)}
 =
 F_{\ell}
 =
 F_{\i}
 =
 \det(1+fL)_{\ell^2(\Z)},
 \end{align*}
 under the assumptions given at the beginning of the section.
 \end{proof}

 Notice that one can significantly relax conditions on parameters set at the beginning of the section for Theorem~\ref{t31}. For this, one could relate matching of determinants to identity \eqref{t12-1}. Since the latter is an equality between polynomials of bounded degree in variables $a_1,\dots,a_N$, $b_1,\dots,b_N$ whose coefficients are rational functions in $q$, it can be extended to all choices of~$a_i$,~$b_j$ and in particular to all choices for which determinants in \eqref{31} make sense.

\subsection[Proof of Lemma 3.4 (i)]{Proof of Lemma~\ref{l32}\,(i)}
\label{s32}
 The case of $\ell=0$, i.e., $F_0=F$, was already mentioned after~\eqref{d31-2}.
 It is then sufficient to show~${F_{\ell}=F_{\ell-1}}$ for any $\ell=1,2,\dots$.
 For this purpose, we rewrite $F_{\ell}$ \eqref{d31-1} as a determinant with rank $N(\ell+ 1)$.
 \begin{Lemma}
 \label{l34}
 We have
 \begin{align}
 \label{3-12}
 F_{\ell}=\det\left(W_{n,n'}\right)_{n,n'=1,2,\dots,N(\ell+1)}.
 \end{align}
 Here $W_{n,n'}$ is defined by
 \begin{align}
\label{3-17}
 W_{n,n'}
 =
 \underset{w=\tilde{a}_{n}^{-1}}{\mathrm{Res}}
 \prod_{i=1}^N\frac{(b_i/w;q)_{\i}}{(a_i w;q)_{\i}}
 \int_{C_{n}}\frac{{\rm d}z}{z}
 \frac{-1}{z-\tilde{a}_{n'}^{-1}}
 \prod_{i=1}^N\frac{(a_i z;q)_{\i}}{(b_i/z;q)_{\i}}
 \frac{\theta\bigl(-\tilde{a}_{n}z/tq^{\frac12+k}\bigr)(q;q)_{\i}^2}{\theta\bigl(-1/tq^{\frac12+k}\bigr)\theta(\tilde{a}_n z)},
\end{align}
 where $C_n$ is a positively oriented circle with the condition $\tilde{a}_n^{-1}<|z|<(q\tilde{a}_n)^{-1}$ and $\theta(x)$ is defined below \eqref{a2-6}.
 \end{Lemma}

 \begin{proof}
 Let us write \eqref{3-10} as $f K_{\ell} = A_{\ell} B_{\ell}$, so that
 \[
 F_{\ell}=\det(1-fK_{\ell})_{\ell^2(\Z)}=\det(1-A_{\ell}B_{\ell})_{\ell^2(\Z)}.
 \]
 By Proposition~\ref{p33}\,(i), $A_{\ell}$, $B_{\ell}$ are
 Hilbert--Schmidt operators and hence $B_{\ell}A_{\ell}$ is trace class.
 Thus we can use the basic property of the determinant \smash{$\det(1-A_{\ell}B_{\ell})_{\ell^2(\Z)}=\det(\delta_{i,j}-(B_{\ell}A_{\ell})_{i,j})_{i,j=1}^{N(\ell+1)}$}
 and rewrite $F_{\ell}$ as
 \begin{align}
 \label{3-17d}
 F_{\ell}=\det(1-A_{\ell}B_{\ell})_{\ell^2(\Z)}
 =\det(\delta_{n,n'}-(B_{\ell}A_{\ell})_{n,n'})_{n,n'=1}^{N(\ell+1)}.
 \end{align}
Here
\begin{align*}
 (B_{\ell}A_{\ell})_{n,n'}
 &=
 \underset{w=\tilde{a}_n^{-1}}{\text{Res}}
 \prod_{i=1}^N\frac{(b_i/w;q)_{\i}}{(a_i w;q)_{\i}}
 \int_C\frac{{\rm d}z}{z}
 \frac{1}{z-\tilde{a}_{n'}^{-1}}
 \prod_{i=1}^N\frac{(a_i z;q)_{\i}}{(b_i/z;q)_{\i}}
 \sum_{m\in\Z}\frac{f(m)}{(z\tilde{a}_{n})^m}
 \\
 &=
 \underset{w=\tilde{a}_{n}^{-1}}{\text{Res}}
 \prod_{i=1}^N\frac{(b_i/w;q)_{\i}}{(a_i w;q)_{\i}}
 \int_{C_{n}'}\frac{{\rm d}z}{z}
 \frac{1}{z-\tilde{a}_{n'}^{-1}}
 \prod_{i=1}^N\frac{(a_i z;q)_{\i}}{(b_i/z;q)_{\i}}
 \cdot
 \frac{\theta\bigl(-\tilde{a}_{n}z/tq^{\frac12+k}\bigr)(q;q)_{\i}^2}{\theta\bigl(-1/tq^{\frac12+k}\bigr)\theta(\tilde{a}_n z)},
\end{align*}
 where in the second equality we applied a version of the Ramanujan summation formula \eqref{a2-6} and $\theta(x)$ is defined below \eqref{a2-6}. Note that to apply
 \eqref{a2-6}, we extend the contour $C$ to $C'_n$,
 which represents the circle centered at the origin with radius $R\in(q/\tilde{a}_n,1/\tilde{a}_n)$.

 Now we show that the entry $\delta_{n,n'}-(BA)_{n,n'}$ is equal to $W_{n,n'}$ \eqref{3-17}. For this, we further extend the radius $R$ of the contour such that $R\in(1/\tilde{a}_n,1/q\tilde{a}_n)$. Recalling the notation \eqref{s32-1} and the definition of $\theta (x)$ introduced below \eqref{a2-6} and
 noting \smash{$
 \prod_{i=1}^N(a_iz;q)_{\infty}=\prod_{\ell=1}^{\infty}(1-\tilde{a}_{\ell}z)$},
 we see
 \begin{gather}
 \underset{z=\tilde{a}_{n}^{-1}}{\text{Res}}
 \frac{{\rm d}z}{z}
 \frac{1}{z-\tilde{a}_{n'}^{-1}}
 \prod_{i=1}^N\frac{(a_i z;q)_{\i}}{(b_i/z;q)_{\i}}
 \cdot
 \frac{\theta\bigl(-\tilde{a}_{n}z/tq^{\frac12+k}\bigr)(q;q)_{\i}^2}{\theta\bigl(-1/tq^{\frac12+k}\bigr)\theta(\tilde{a}_n z)}\nonumber \\
 \qquad=\underset{z=\tilde{a}_{n}^{-1}}{\text{Res}}
 \frac{-\tilde{a}_{n'}}{z}
 \frac{\prod_{\substack{\ell=1\\ \ell\neq n'}}^{\infty}(1-\tilde{a}_\ell z)}{\prod_{i=1}^N(b_i/z;q)_{\i}}
 \cdot
 \frac{\theta\bigl(-\tilde{a}_{n}z/tq^{\frac12+k}\bigr)(q;q)_{\i}^2}{-\tilde{a}_n\bigl(z-\tilde{a}_n^{-1}\bigr)(q\tilde{a}_nz;q)_\infty \bigl(q\tilde{a}_n^{-1}/z;q\bigr)_\infty\theta\bigl(-1/tq^{\frac12+k}\bigr)}
 \nonumber\\
 \qquad = \tilde{a}_n
 \frac{\prod_{\substack{\ell=1\\ \ell\neq n'}}^{\infty}\bigl(1-\tilde{a}_\ell \tilde{a}_n^{-1}\bigr)}{\prod_{i=1}^N(b_i \tilde{a}_n;q)_{\i}}\delta_{n,n'}. \label{eq:zres}
 \end{gather}
 Here in the second equality we used the fact that in the second expression all factors do not include the residue at $z=\tilde{a}_n^{-1}$ except the factor $1/\big(z-\tilde{a}_n^{-1}\big)$ and this is eliminated by the factor~\smash{${\prod_{\substack{\ell=1\\ \ell\neq n'}}^{\infty}(1-\tilde{a}_\ell z)}$} unless $n=n'$.
 Similarly, we compute
 \begin{align}\label{eq:wres}
 \underset{w=\tilde{a}_{n}^{-1}}{\text{Res}}
 \prod_{i=1}^N\frac{(b_i/w;q)_{\i}}{(a_i w;q)_{\i}}=\underset{w=\tilde{a}_{n}^{-1}}{\text{Res}}\frac{\prod_{i=1}^N(b_i/w;q)_{\i}}{\prod_{\ell=1}^{\infty}(1-\tilde{a}_{\ell}w)}=-\tilde{a}_n^{-1}
 \frac{\prod_{i=1}^N(b_i \tilde{a}_n;q)_{\i}}{\prod_{\substack{\ell=1\\ \ell\neq n}}^{\infty}\big(1-\tilde{a}_\ell \tilde{a}_n^{-1}\big)}.
 \end{align}
 Therefore from \eqref{eq:zres} and \eqref{eq:wres}, the extension of the contour cancels the term $\delta_{n,n'}$ in the last expression in \eqref{3-17d} and we arrive at the expression \eqref{3-17}.
 \end{proof}

 Now we prove Lemma~\ref{l32}\,(i). We focus on the representation \eqref{3-12} and especially on the last ($N(\ell+1)$-th) row of the matrix $W$. We decompose the elements into two parts as
\begin{align*}
 W_{N(\ell+1),n'}
 &=W^{(1)}_{N(\ell+1),n'}+W^{(2)}_{N(\ell +1),n'}
\end{align*}
 for $n'\!=\!1,\dots,N(\ell+1)$. Here
 \smash{$W^{(1)}_{N(\ell+1),n'}$} corresponds to
 the residue at \smash{$z\!=\!\tilde{a}_{N(\ell+1)}^{-1}$}
 while~\smash{$W^{(2)}_{N(\ell+1),n'}$} corresponds to the element with the shrunk contour
 satisfying
 $q\tilde{a}_{n}^{-1}<|z|<\tilde{a}_{n}^{-1}$.
 As we stated in the last paragraph in the proof of Lemma~\ref{l34}, one easily sees \smash{$W^{(1)}_{N(\ell+1),n'}=\delta_{N(\ell+1),n'}$}.
 On the other hand recalling the notation \eqref{s32-1} and noting $q^{-1}\tilde{a}_n=\tilde{a}_{n-N}$, we have
 \begin{gather}
\underset{w=\tilde{a}_{n}^{-1}}{\text{Res}}
 \prod_{i=1}^N\frac{(b_i/w;q)_{\i}}{(a_i w;q)_{\i}}
 = q^{-1}
 \prod_{i=1}^N\frac{1}{1-b_i\tilde{a}_nq^{-1}}\frac{1}{1-a_i\tilde{a}_n^{-1}}
 \cdot
 \underset{w=\tilde{a}_{n-N}^{-1}}{\text{Res}}
 \prod_{i=1}^N\frac{(b_i/w;q)_{\i}}{(a_i w;q)_{\i}}, \label{eq:wn1}
 \\
 \theta(\tilde{a}_n x)=\frac{-q}{\tilde{a}_n x}
 \theta(\tilde{a}_{n-N}x).\label{eq:wn2}
 \end{gather}
 Applying \eqref{eq:wn1} to the first factor in \eqref{3-17} and \eqref{eq:wn2} to two factors $\theta\bigl(-\tilde{a}_n z/tq^{\frac12+k}\bigr)$ and $\theta(\tilde{a}_nz)^{-1}$ in the same equation, we see
 \[
 W^{(2)}_{N(\ell+1),n'}=-tq^{-\frac12+k}\prod_{i=1}^N\frac{1}{1-b_i\tilde{a}_{N(\ell+1)}q^{-1}}\frac{1}{1-a_i\tilde{a}_{N(\ell+1)}^{-1}}\cdot W_{N\ell,n'}.
 \]
 Combining them
 with the multi-linearity of determinants,
 we get
\[
 F_{\ell}=\det[W_{n,n'}]_{n,n'=1}^{N(\ell +1)}
 =
 \det[W_{n,n'}]_{n,n'=1}^{N(\ell+1)-1}.
\]
 Repeating this procedure $N-1$ times, we arrive at $F_{\ell}=F_{\ell-1}$.

\subsection[Proof of Lemma 3.4 (ii)]{Proof of Lemma~\ref{l32}\,(ii)}
\label{s33}
 From \eqref{3-10}
 and Corollary~\ref{cb2}, we find a uniform bound for the kernel
 \[
 \t(m_1)f(m_1)K_{\ell}(m_1,m_2)\t^{-1}(m_2),
 \]
 where
 $K_{\ell}(m_1,m_2)$, $\t(m)$ and $f(m)$ are defined by \eqref{d31-3}, \eqref{3-12d} and \eqref{p251-1} respectively. We have
 \begin{align}
 \label{3-32}
 \big|\t(m_1)f(m_1)K_{\ell}(m_1,m_2)\t^{-1}(m_2) \big|\le Dd^{|m_1|+|m_2|}
 \end{align}
 for some constants $D>0$ and $d\in (0,1)$ which do {\it not} depend on $\ell$, $m_1$ or $m_2$.

 Lemma~\ref{l32}\,(ii) is then immediately proven
 combining \eqref{3-32} with
 the Hadamard inequality
\begin{align}
\label{33-1}
 |\det A|\le
 \prod_{i=1}^n\Bigg(\sum_{j=1}^n|a_{i,j}|^2\Bigg)^{1/2},
\end{align}
which holds for any for $n\times n$ matrix $A=(a_{ij})_{i,j=1,\dots,n}$.
 Note that $F_\ell$ \eqref{d31-1} can be written as
\begin{align}
 F_{\ell}
\label{33-2}
 &=\sum_{n=0}^{\i}\frac{(-1)^n}{n!}
 \sum_{m_1=-\i}^{\i}\cdots\sum_{m_n=-\i}^{\i}
 \det\big(\t(m_i)\t^{-1}(m_j)f(m_i)K_\ell (m_i,m_j)\big)_{i,j=1}^n.
\end{align}
 By \eqref{3-32} and \eqref{33-1},
 the determinant in \eqref{33-2} can be bounded as
\begin{gather}
 \big|\det\big(\t(m_i)\t^{-1}(m_j)f(m_i)K_\ell (m_i,m_j)\big)_{i,j=1}^n \big|
 \nonumber\\
 \qquad \le
 \prod_{i=1}^n\Bigg(\sum_{j=1}^n \big|\t(m_i)\t^{-1}(m_j)f(m_i)K_\ell (m_i,m_j)\big|^2
 \Bigg)^{1/2}
 \nonumber\\
 \qquad\le
 \prod_{i=1}^n\Bigg(\sum_{j=1}^n \big|Dd^{|m_i|+|m_j|}\big|^2
 \Bigg)^{1/2}
 =D^n \prod_{i=1}^n d^{|m_i|}\Bigg(\sum_{j=1}^n d^{2|m_j|}\Bigg)^{1/2}\nonumber\\
 \qquad \le D^n n^{\frac{n}{2}}\prod_{i=1}^n d^{|m_i|}.\label{33-3}
\end{gather}
 Note that bound \eqref{33-3} holds for any $\ell\in\Z_{\ge 0}$. Moreover, we have
\begin{gather*}
 \sum_{m_1=-\i}^{\i}\cdots\sum_{m_n=-\i}^{\i}
 D^n n^{\frac{n}{2}}\prod_{i=1}^n d^{|m_i|}
 =D'^n n^{\frac{n}{2}} <\i,
 \qquad \sum_{n=0}^{\i}\frac{1}{n!}D'^n n^{\frac{n}{2}}<\i
\end{gather*}
with some constant $D'>0$.
 Thus by dominated convergence theorem and Proposition~\ref{p33}\,(ii), Lemma \ref{l32}
 (ii) holds.
\qed

\subsection[Proof of Lemma 3.4 (iii)]{Proof of Lemma~\ref{l32}\,(iii)}
\label{s34}

\begin{figure}[ht]
 \centering
 \includegraphics[scale=1]{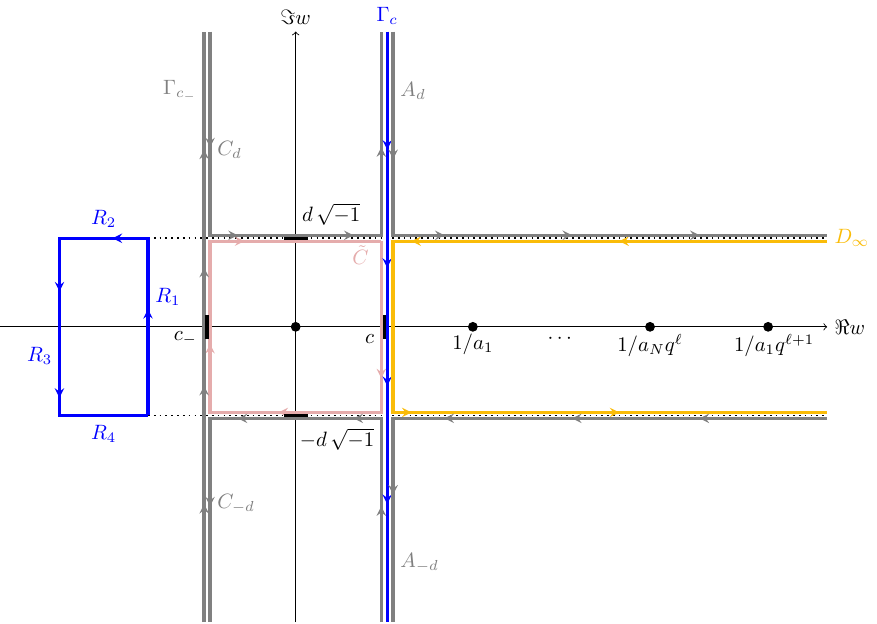}
 \caption{All of the contours appearing in the proof of
 Lemma \ref{l32}. We set
 $c>0$, $c_-<0$. }
 \label{f3}
\end{figure}

 In this subsection, we will
 prove $K_{\i}(m_1,m_2)
 =-L(m_1,m_2)$ by deforming the contour of $w$
 from $D_{\infty}$ to $\tilde{C}$ going through
 $\Gamma_c$, where all contours are depicted in Figure~\ref{f3}. We have to care about the behavior of the integrand $g_{a,b}(z,w;m_1,m_2)$ \eqref{p223}
 when $|w|$ is large since the poles~${w=1/a_iq^j}$,~${i=1,\dots,N}$, $j\in\Z_{\ge 0}$ accumulate at the infinite point of $w$ producing an essential singularity.
 We estimate the behavior for large $|w|$ using Appendix~\ref{ac}.

 Here for convenience we also use
 the notation \eqref{3-26} where
 we set $\sharp$ to be some contours
 depicted in
 Figure~\ref{f3}. Note that by definition, one sees
 $K_{\infty}(m_1,m_2)= K(m_1,m_2;D_{\i})$,
 and ${L(m_1,m_2)=-K\big(m_1,m_2;\tilde{C}\big)}$.
 Note that $\tilde{C}$ goes around the origin in negative direction.
 We now show that,
 from the properties of the $q$-Pochhammer symbol discussed in
 Appendix \ref{ac}, the following holds,
\begin{align}
\label{s34-2}
 K(m_1,m_2;\sharp)=0,
 \qquad \text{when}\quad \sharp=A_d,A_{-d},C_d,C_{-d}\ \text{and}\ \Gamma_{c_-}.
\end{align}
 In Figure~\ref{f3}, these contours
 are displayed in grey.
 These facts can be
 seen in the following way.
 In \eqref{3-26}, we focus on the
 integration of $w$,
 \begin{align}
 \label{3-29}
 \int_{\sharp}{\rm d}w\frac{w^{m}}{z-w}
 \prod_{i=1}^N
 \frac{{(b_i/w;q)_{\i}}}{(a_iw;q)_{\i}},
 \end{align}
 with $m\in\Z$ and $z\in C$ fixed. Here the contour $C$ is illustrated in Figure~\ref{f1}.
 Hereafter, we will explain the case $\sharp=\Gamma_{c_-}$ in detail whereas in other cases we will give an outline only in the paragraph below \eqref{3-34} since
 the other cases follow in a similar
 way.

 First we note that \eqref{3-29} does not depend on the value of $c_-$, as can be seen as follows. We consider the contour integral \eqref{3-29}
 with $\sharp=R_1\cup R_2\cup R_3\cup R_4$ where $R_i$, $i=1,2,3,4$ are depicted in
 Figure~\ref{f3}. Since there are no singularities inside the rectangle, the value of this integral is zero by the Cauchy's theorem. On the other hand
 from Lemma~\ref{lc1}\,(ii) we find the contributions from both $R_2$ and $R_4$ are bounded by $cd^m \exp\bigl[-Nc'\log^2 d\bigr]$ when $d>1$ for some constants $c$ and~$c'$, thus they vanish when $d\rightarrow\i$. Combining these facts, we find
 the contribution from $R_1$ is equal to that from $R_3$
 as $d\rightarrow\i$ up to sign.

 Next, we will show that \eqref{3-29} with $\sharp=\Gamma_{c_-}$ vanishes when $c_-\rightarrow -\i$. Combining Lem\-ma~\ref{lc1}\,(i) and (ii), we see that
 for $z=a+b\sqrt{-1}$ with $a<-1$ and $b\in\R$,
 \begin{align}
 \label{3-33}
 |(z;q)_{\i}|\ge c_1\exp\bigl[c_2\log^2(|a|\vee |b|)\bigr],
 \end{align}
 where $x\vee y=\max(x,y)$. Applying \eqref{3-33} to
 \eqref{3-29} with $\sharp=\Gamma_{c_-}$, we have
 \begin{align}
 \bigg| \int_{\Gamma_{c_-}} {\rm d}w h(w)\bigg|
 \le{}& \int_{\Gamma^{(1)}_{c_-}}{\rm d}w |h(w)|+\int_{\Gamma^{(2)}_{c_-}}{\rm d}w |h(w)|
 \nonumber\\
\le{}& c_3 |c_-|^{m_2}\exp\bigl(-Nc_2\log^2(|c_-|)\bigr) \nonumber\\
 &+c_4\int_{\Gamma_{c-}^{(2)}} {\rm d}w |w|^{m_2}
 \exp\bigl(-Nc_2\log^2|w|\bigr),
 \label{3-34}
 \end{align}
 where $h(w)$ is the integrand of \eqref{3-29}, $\Gamma_{c_-}^{(1)}$
 is the vertical line from $c_-+c_-\sqrt{-1}$ to $c_--c_-\sqrt{-1}$ and \smash{$\Gamma_{c_-}^{(2)}$} is its complement \smash{$\Gamma_{c_-}\setminus\Gamma_{c_-}^{(1)}$}.
 Thus we find that in the last expression in \eqref{3-34}, both two terms vanishes as $|c_-|\rightarrow\i$ and we arrive at the conclusion.

 We can also obtain \eqref{s34-2} with the other contours $\sharp=A_d,A_{-d},C_d,C_{-d}$ in a similar way, i.e., they follow from the following two facts: (I) The integration value \eqref{3-29} does not depend on the value of $d>1$ for any $m\in\mathbb{Z}$ and $z\in C$. (II) \eqref{3-29} vanishes as $d\rightarrow\infty$ for any~${m\in\mathbb{Z}}$ and $z\in C$. For proving (I), as in the case of $\Gamma_{c-}$, we introduce a positively-oriented rectangular contour~${\tilde{R}:=\tilde{R}_1\cup\tilde{R}_2\cup\tilde{R}_3\cup\tilde{R}_4}$, where $\tilde{R}_1$, $\tilde{R}_2$, $\tilde{R}_3$ and $\tilde{R}_4$ represents the right, top, left, and bottom edge respectively. In the case $A_d$ (resp.\ $A_{-d}$), we put $\tilde{R}$ such that its lower-left corner match the corner of $A_d$ (resp.\ $A_{-d}$). Then by using \eqref{lc1-3} in Lemma~\ref{lc1}\,(iii), we see that \eqref{3-29} with $\sharp=\tilde{R}_1$ vanishes as the length of the horizontal sides goes to infinity. In the case $C_d$ (resp.~$C_{-d}$), we set the length of the horizontal sides in $\tilde{R}$ to be $c-c_-$ and place it in a way that its two bottom (resp.\ top) corners match the corners in $C_d$ (resp.\ $C_{-d}$) and show that~\eqref{3-29} with~${\sharp=\tilde{R}_2}$ (resp.\ $\sharp{R}_4$) vanishes as the length of the vertical sides goes to infinity using Lemma~\ref{lc1} (ii). On the other hand one can show the property (II) just by using Lemma~\ref{lc1}\,(ii) and \eqref{lc1-3} in Lemma~\ref{lc1}\,(iii) in the cases $A_d$ and $A_{-d}$ while in the case $C_{d}$ and~$C_{-d}$, we use Lemma~\ref{lc1}\,(ii).

 Now we will prove Lemma \ref{l32}\,(iii)
 by showing the following two relations:
\begin{subequations}
 \begin{align}
\label{s34-5}
 &K(m_1,m_2;\Gamma_c)
 =K_{\i}(m_1,m_2),
 \\
\label{s34-6}
 &K(m_1,m_2;\Gamma_c)
 =-L(m_1,m_2).
 \end{align}
\end{subequations}
 First, we prove \eqref{s34-5}. Noting $\Gamma_c=A_d\cup D_{\infty}\cup A_{-d}$, and \eqref{s34-2}, we have
\begin{align*}
 K(m_1,m_2;\Gamma_c)
 &=K(m_1,m_2;D_{\infty})+K(m_1,m_2;A_d)+K(m_1,m_2;A_{-d})\\
 &=K(m_1,m_2;D_{\infty}).
\end{align*}
 Recalling by definition $K(m_1,m_2; D_{\infty})=K_{\infty}(m_1,m_2)$, we arrive at \eqref{s34-5}.

 Next we show the relation \eqref{s34-6}.
 Noting the decomposition
 $\tilde{C}=\Gamma_{c}\cup C_d\cup C_{-d}\cup\Gamma_{c_-}$, and~\eqref{s34-2}, we have
 \[
 K\big(m_1,m_2,\tilde{C}\big)=
 K(m_1,m_2,\Gamma_c).
 \]
 Combining this with the fact $L(m_1,m_2)=-K\big(m_1,m_2;\tilde{C}\big)$, we arrive at \eqref{s34-6}.
\qed

\section{Comparison with Borodin's matching} \label{sec:comparison}

This section stems from conversations we had in separate occasions with Alexei Borodin and Guillaume Barraquand, soon after posting the first version of this paper on arXiv. We are going to compare results proven in the previous sections and especially Theorem~\ref{thm:intro_prob} with the matching between the distribution of the length of a random partition in the Hall--Littlewood measure and a certain multiplicative statistics of the Schur measure proven in \cite{Borodin2018momentmatch}. The Hall--Littlewood measure is defined by
\begin{align}\label{eq:HLmeasure}
 M_{\text{HL}}(\mu)&=P_{\mu}(a_1,\dots,a_N;0,t)Q_{\mu}(b_1,\dots,b_M;0,t)/Z_{\text{HL}},
\end{align}
where $P_{\mu}(a_1,\dots,a_N;0,t)$ and $Q_{\mu}(b_1,\dots,b_M;0,t)$ are the Hall--Littlewood polynomials, which are the case $q=0$ of the Macdonald polynomials $P_{\mu}(a_1,\dots,a_N;q,t)$ and $Q_{\mu}(b_1,\dots,b_M;q,t)$ defined in \cite[Chapter VI] {macdonald1998symmetric}
 and $Z_{\text{HL}}=\prod_{i=1}^N\prod_{j=1}^M (1-ta_ib_j)/(1-a_ib_j)$ is the normalization constant.

\begin{Proposition}[\cite{Borodin2018momentmatch}]
Let $a$, $b$ be specializations of the algebra of symmetric functions such that $0<|Z_{\mathrm{HL}}|<\infty$ and $0<|Z_{\mathrm{S}}|<\infty$, where \[Z_{\mathrm{HL}}=\sum_{\mu}P_{\mu}(a;0,t)Q_{\mu}(b;0,t)\qquad \text{and}\qquad Z_{\mathrm{S}}\!=\!\sum_{\lambda}s_{\lambda}(a)s_{\lambda}(b).
\]
 Let $HL(a,b)$ and $S(a,b)$ denote respectively the Hall--Littlewood and the Schur $($possibly signed$)$ measure.
Let $\mathbb{E}_{HL(a,b)}[\cdot]$ and $\mathbb{E}_{S(a,b)}[\cdot]$ denote
the average over the Hall--Littlewood measure $M_{\text{HL}}(\mu)$ \eqref{eq:HLmeasure} and the Schur measure $M_{\text{S}}(\mu)$ \eqref{1-15} respectively.
Then we have
	 \begin{equation} \label{matching_hl_schur}
 	\mathbb{E}_{HL(a,b)}\left[ \frac{1}{\bigl(-vt^{-\ell(\mu)};t\bigr)_\infty} \right] = \mathbb{E}_{S(a,b)}\Bigg[ \prod_{i\ge 1} \frac{1}{1+v t^{-1+i-\mu'_i}} \Bigg].
 \end{equation}
\end{Proposition}

In his original paper \cite[Corollary 4.4]{Borodin2018momentmatch}, Borodin proves a variant of identity \eqref{matching_hl_schur} relating certain averages with respect to the stochastic six vertex model and the Macdonald measure. One can replace the average with respect to the stochastic six vertex model with the Hall--Littlewood measure thanks to the equivalence in law between the height function and the length of a Hall--Littlewood random partition (this equivalence was proven in \cite{BorodinBufetovWheeler2016}). On the other hand one can also replace the average over the Macdonald measure with the Schur measure since the identity holds for arbitrary $t$, which is the one of the parameters $(q,t)$ in Macdonald polynomial thus one can set $t=q$. From these properties, we see that the identity in \cite[Corollary 4.4]{Borodin2018momentmatch} leads to \eqref{matching_hl_schur}. Below we report a different proof that was suggested to us by Guillaume Barraquand and that to the best of our knowledge has not appeared in published literature before, although it has been known for a few years.\footnote{For instance, the same argument we present here had been used in 2019 by Corwin in the lecture series in Kyushu \url{https://www2.math.kyushu-u.ac.jp/~osada-labo/msj-si2019/slides/Corwin2.pdf}.}

\begin{proof} Hereafter, we set $N=M=n$ in \eqref{eq:HLmeasure} and \eqref{1-15}.
	It was observed in \cite{macdonald1998symmetric} that the summation of Macdonald polynomials
 \begin{equation} \label{eq:q-independent expression}
 \frac{1}{\Pi(a,b;q,t)} \sum_\mu \prod_{i=1}^n\big(1-uq^{\mu_i} t^{n-i}\big) P_\mu(a;q,t) Q_\mu(b;q,t),
 \end{equation}
 where
 \[
 \Pi(a,b;q,t) = \sum_{\mu} P_\mu(a;q,t) Q_\mu(b;q,t) = \prod_{i,j=1}^n \frac{(t a_i b_j;q)_\infty}{(a_i b_j;q)_\infty}
 \]
 does not depend on $q$, when $a=(a_1,\dots,a_n)$ and $b=(b_1,\dots,b_n)$ are $n$-tuples of complex numbers. More precisely, the proof of the $q$-independence of \eqref{eq:q-independent expression} is contained in the proof of \cite[Chapter VI, Section~3, equation~(3.12)]{macdonald1998symmetric}. Setting $q=t$, this fact produces a matching between expectations of Schur and Macdonald measures (the latter denoted by $M(a,b)$)
 \[
 \mathbb{E}_{M(a,b)}\left[ \prod_{i=1}^n \big(1-uq^{\mu_i} t^{n-i}\big) \right] = \mathbb{E}_{S(a,b)}\left[ \prod_{i=1}^n \big(1-u t^{\mu_i + n-i}\big) \right].
 \]
 Setting $q=0$, the left-hand side becomes an expectation over the Hall--Littlewood measure and after simple algebraic manipulation, we obtain
 \begin{equation*}
 \mathbb{E}_{HL(a,b)}\Bigg[ \frac{1}{\big(ut^{n-\ell(\mu)};t\big)_\infty} \Bigg] = \mathbb{E}_{S(a,b)}\Bigg[ \frac{ \prod_{i=1}^n \big(1 - u t^{\mu_{i} - i + n}\big) }{ \prod_{i\ge0} \big(1- u t^{i} \big) } \Bigg]
 =
 \mathbb{E}_{S(a,b)}\Bigg[ \prod_{i\ge 1} \frac{1}{1-u t^{n-1+i-\mu'_i}} \Bigg].
 \end{equation*}
 We can finally rescale the $u$ parameter as $u=-vt^{-n}$ to obtain the equality \eqref{matching_hl_schur}, which holds for $a$, $b$, being $n$-tuples of complex numbers of modulus smaller than 1. Note that \eqref{matching_hl_schur} can be regarded as an identity for the symmetric polynomials with $n$ variables. In addition using stability properties (i.e., the property that the identity for variables $(a_1,\dots,a_{n-1})$ is equal to the identity for $(a_1,\dots,a_{n-1},a_n=0)$) 
 of the Hall--Littlewood and Schur polynomials we can extend the identity of the symmetric polynomials to that of the
 symmetric functions. Thus we see that \eqref{matching_hl_schur} holds in the algebra of symmetric functions; for more information about stability and inverse limit procedures, see \cite[Chapter~I.2]{macdonald1998symmetric}. Through this extension $a$ and $b$ become specializations of the algebra of symmetric functions and we assume the condition that the Macdonald measure is absolutely summable. This completes the proof.
\end{proof}

We would now like to transform the identity \eqref{matching_hl_schur} into an analogous identity relating the $q$-Whittaker measure and the Schur measure. For this, we introduce the specialization of the algebra of symmetric functions defined by its action on (the algebraic basis of) power sum symmetric functions $p_n$ as
	\[
 		\widetilde{a}\colon\ p_n \mapsto \frac{(-1)^{n-1}}{1-t^n} p_n(a),
	\]
	where $a$ represents a collection of complex numbers and we recall that, under such specialization the Hall--Littlewood functions are turned into $q$-Whittaker functions (with $q=t$) as $P_\mu(\widetilde{a};0,t) = P_{\mu'}(a;t,0)$.
	Considering the Hall--Littlewood measure in \eqref{matching_hl_schur} with specializations set to $\widetilde{a}$, \smash{$\widetilde{b}$}, we get\looseness=-1
	\begin{equation} \label{matching_qW_schur}
		\mathbb{E}_{t\mathrm{W}(a,b)}\left[ \frac{1}{(-vt^{-\mu_1};t)_\infty} \right] = \mathbb{E}_{S(a',b')}\Bigg[ \prod_{i\ge 1} \frac{1}{1+v t^{-1+i-\mu_i}} \Bigg],
	\end{equation}
	where $t\mathrm{W}(a,b)$ denotes the $t$-Whittaker measure. Above the specialization of the Schur measure~is
	\[
		a' \colon\ p_n \mapsto \frac{1}{1-t^n} p_n(a)
	\]
	and we used the fact that
\[
s_\mu(\widetilde{a}) = s_{\mu'}(a')=s_{\mu'}(\alpha),
\]
where in the right-hand side $\alpha$
is the collection of the complex numbers whose set $\{\alpha_i\}_{i=1,2,\dots}$ is given by
\[
 \{\alpha_i\}_{i=1,2,\dots}=\bigsqcup_{j=1,2,\dots}\big\{a_j,a_jt,a_jt^2,\dots\big\}.
\]

Through \eqref{matching_qW_schur} one can write Fredholm determinant expressions for the $t$-Laplace transform of the probability distribution of the first row of a partition $\mu$ distributed according to the $t$-Whittaker measure. This is because the right-hand side of $\eqref{matching_qW_schur}$ is a multiplicative observable of a determinantal point process, namely $(\lambda_i-i+1/2)_{i=1,2,\dots}$ when $\lambda$ is distributed according to the Schur measure. We have, changing all $t$ parameters to $q$
 \begin{equation} \label{eq:fredholm_det_qW_schur}
 \mathbb{E}_{q\mathrm{W}(a,b)}\left[ \frac{1}{(-vq^{-\mu_1};q)_\infty} \right] = \det \big( 1- \bar{f}_v L^* \big)_{\ell^2(\mathbb{Z})},
 \end{equation}
 where $\bar{f}_v(x) = \frac{ 1 }{1 + v q^{x+1/2}}$
 and $L^*(i,j)=-L(-1-j,-1-i)$, where $L$ was defined in \eqref{t252}. To check that $L^*$ is the correlation kernel of a Schur measure (with specializations $a'$, $b'$) see, for instance \cite[Theorem~1.1]{aggarwal2015correlation_schur}. Now, simple algebraic manipulations of the right-hand side of \eqref{eq:fredholm_det_qW_schur} show that
 \begin{equation} \label{eq:det_L_star_L}
 	\det \big( 1- \bar{f}_v L^* \big)_{\ell^2(\mathbb{Z})} = \det ( 1+ f_{-vq^{1/2}} L )_{\ell^2(\mathbb{Z})},
 \end{equation}
 where the function $f_\zeta$ was defined in Proposition~\ref{p22}. We can now compare the expression in the right-hand side of \eqref{eq:det_L_star_L} with that in the right-hand side of \eqref{p251} to discover that they differ by a change of parameter $v\to t^k$.

The argument outlined in the previous paragraphs yields another proof of the correspondence between the first row of the $q$-Whittaker measure and the edge of the periodic Schur measure.
Furthermore, we can also obtain \eqref{matching_hl_schur} from our relation \eqref{t1-1}. Noting \eqref{t1-1} holds for general specialization of symmetric functions as stated in Remark~\ref{rem:intro} and applying the specialization~${p_n\mapsto (-1)^{n-1}(1-t^n)p_n(a)}$ to our relation, we arrive at \eqref{matching_hl_schur}. This shows the equivalence between Borodin's matching \eqref{matching_hl_schur} and our Theorem~\ref{thm:intro_prob}.

\appendix
\section{Notations and formulas}
\label{aa}
 In this sections, we summarize the notations
 and formulas which we frequently use in this paper.
\subsection{Partitions}
\label{a1}
 A partition $\l$
 is a sequence of nonnegative integers
 $\l=(\l_1,\l_2,\dots)$, $\l_j\in\Z_{\ge 0}$ with non decreasing orders $\l_1\ge\l_2\ge\cdots$. The length
 of partition $\ell(\l)$ is the number of nonzero
 elements in $\l$ and the size of $\l$ is defined
 as $|\l|=\sum_{j}\l_j$. Let $\Pe$ be the set of
 the partitions whose lengths are finite.
 The partition $\l\in\Pe$
 can be expressed in another way
 $\l=1^{m_1}2^{m_2}\cdots$, where~$m_j$ represents the number of
 the elements $\l_k$s such that
 $\l_k=j$. Note that $\sum_jm_j=n$
 and $|\l|=\sum_{j} jm_j$.
 For example, $\l=(3,3,2,2,1,1,1,0,\dots)$ is
 expressed as $1^3 2^2 3^2$. We write $\lambda=\varnothing$ if all $\lambda_j$'s are~0's.
 For $\l,\mu\in\Pe$, we denote $\mu\subset\l$ if
 $\mu_i\le\l_i,~i=1,2,\dots$.

\subsection[q-Series and sum over partitions]{$\boldsymbol{q}$-Series and sum over partitions}
\label{a2}
 Here we summarize a few formulas for $q$-series and related summations over partitions. For more detail, we refer the reader to \cite[Sections~10 and~11]{AndrewsAskeyRoy2000}. Define
\begin{align}
\label{a2-1}
 &(a;q)_\i = \prod_{j=0}^\i \big(1-a q^j\big),
 \\
\label{a2-1d}
 &(a;q)_n = \frac{(a;q)_\i}{(a q^n;q)_\i}
 =
 \begin{cases}
 (1-a)(1-aq)\cdots \big(1-a q^{n-1}\big), & n\ge 0,
 \\
 \dfrac{1}{(1-a/q)(1-a/q^2)\cdots(1-a/q^m)},& n=-m,~m>0.
 \end{cases}
\end{align}
Note that $(q;q)_{-m}^{-1}=0$ for $m\in\Z_{>0}$.

 Next we give some formulas for the $q$-series. First, we state the $q$-binomial theorem:
 For $|z|<1$, $|q|<1$ and $a\in\C$,
\begin{align}
\label{a2-2}
\sum_{n=0}^\i \frac{(a;q)_n}{(q;q)_n} z^n
=
\frac{(az;q)_\i}{(z;q)_{\i}}.
\end{align}
For the proof, see, e.g., \cite[Theorem 10.2.1]{AndrewsAskeyRoy2000}. In this paper,
we use the two special cases,
\begin{align}
\label{a2-3}
&\sum_{n=0}^\i \frac{z^n}{(q;q)_n}
=
\frac{1}{(z;q)_{\i}},
\\
 \label{a2-4}
&\sum_{n=0}^\i \frac{(-1)^n q^{\frac{n(n-1)}{2}}}{(q;q)_n} z^n
= (z;q)_\i.
\end{align}
 The first equality \eqref{a2-3} holds for $|z|<1$ while the second one \eqref{a2-4} holds for arbitrary $z\in\C$. They are written in~\cite[Corollary
 10.2.2\,(a) and (b)]{AndrewsAskeyRoy2000}.
 One immediately obtains \eqref{a2-3} (resp.~\eqref{a2-4})
 by setting $a=0$ (resp.\ replacing
 $a\rightarrow 1/a$ $x\rightarrow ax$, then taking the limit $a\rightarrow 0$) in \eqref{a2-2}.

 The Ramanujan's summation formula is a
 generalization of the $q$-Binomial
 theorem~\eqref{a2-2} to the sum over
 $\Z$: For $|q|<1$, $|b/a|<|z|<1$,
\begin{equation}
\label{a2-5}
\sum_{n\in\Z} \frac{(a;q)_n}{(b;q)_n} z^n
=
\frac{(az;q)_\i \big(\frac{q}{az};q\big)_\i (q;q)_\i \big(\frac{b}{a};q\big)_\i}
 {(z;q)_\i \big(\frac{q}{a};q\big)_\i (b;q)_\i \big(\frac{b}{az};q\big)_\i}.
\end{equation}
 For a proof, see, e.g., in \cite[Theorem 10.5.1]{AndrewsAskeyRoy2000}.
 In this paper, we use the special case
 $a=-t$, $b=-qt$, $z=q/w$ with
 $t\in \C\setminus\{0\}$,
\begin{align}
 \label{a2-6}
 \sum_{n\in\Z}\frac{t q^n w^{-n}}{1+tq^n}
 =
 \frac{\theta(-w/t)(q;q)_{\i}^2}
 {\theta(-1/t)\theta(w)},
\end{align}
 where $\theta (x)$ is defined by using
 \eqref{a2-1} as $\theta(x)=(x;q)_{\i}(q/x;q)_{\i}$
 and
 according to the condition on $z$ in \eqref{a2-5}, $w$ should satisfy $|q|<|w|<1$.

Finally, we present a well-known formula for the generating function of $\lambda$:
\begin{equation}
 \label{a2-8}
 \sum_{\l\in\Pe}
 q^{|\l|}=\frac{1}{(q;q)_\i}
\end{equation}
and its refinement with a fixed $\lambda_1$:
\begin{align}
\label{a2-7}
 &\sum_{\substack{\l\in\Pe\\ \l_1=n}}
 q^{|\l|}=\frac{q^n}{(q;q)_n}.
\end{align}
 \eqref{a2-7} is obtained as follows. Switching to the notation of partition as $\l=1^{m_1}2^{m_2}\cdots$,
 and noting $|\l|=\sum_{j}jm_j$, we have
\begin{align*}
 \sum_{\substack{\l\in\Pe\\ \l_1=n}}
 q^{|\l|}
 =
 \prod_{j=1}^{n-1}\sum_{m_j=0}^{\i}
 q^{jm_j}\cdot \sum_{m_n=1}^{\i}q^{n m_n}=\prod_{j=1}^{n-1}\frac{1}{1-q^j}
 \cdot\frac{q^n}{1-q^n}.
\end{align*}
\eqref{a2-8} follows immediately from \eqref{a2-7} by taking the sum over $n$ and noting \eqref{a2-3}.

\section[Properties of the kernels A(m;r) and B(r;m)]{Properties of the kernels $\boldsymbol{A(m;r)}$ and $\boldsymbol{B(r;m)}$}
\label{ab}
 In this appendix, we give upper bounds for matrix elements of $A(m;r)$ \eqref{3-15d} and $B(r;m)$ \eqref{3-16d}, which are useful to prove the well-definedness of the Fredholm determinants $F_{\ell}$\eqref{d31-1} and $F_{\i}$\eqref{d31-2} in Proposition~\ref{p33}\,(ii) and Lemma~\ref{l32}\,(ii) in Section~\ref{s33}.
\begin{Lemma}
\label{lb1}
Recall the assumption on $a_i$, $b_i$, $i=1,\dots,N$ stated at the beginning of Section~{\rm\ref{s3}}. Recall also the notation \eqref{s32-1}, $\tilde{a}_r=a_k q^{u}$. Then there exist constants $C_+$, $C_-$, $D_+$, and $D_-$ which do not depend on $m$ or
$r$ such that the following inequalities hold.
\begin{itemize}\itemsep=0pt
 \item [$(i)$] Assume $b$ satisfies $b_{\max}<b<1$. We have
 \begin{align}
 \label{bb2-1}
 |A(m;r)|
 \begin{cases}
 < C_+ a_{1}^mq^{\frac{m^2}{2N}+\frac{m}{2}+u}, & m\ge 0,
 \\
 < C_- b^{-m} q^{u}, & m<0.
 \end{cases}
 \end{align}
 \item[$(ii)$] Assume $a_{1}/a_{N}<q^{-\frac12+\e}$ with $\e\in(0,1/2)$. We have
 \begin{align}
 \label{bb2-2}
 |B(r;m)|
 \begin{cases}
 < D_+ a_{N}^{-m}q^{-\frac{m^2}{2N}+\e m-(\frac12+\e)u}, & m\ge 0,
 \\
 < D_- a_{1}^{-m}q^{u}, & m<0.
 \end{cases}
 \end{align}
\end{itemize}
\end{Lemma}

From \eqref{bb2-1}, \eqref{bb2-2}, we immediately obtain the bounds
for $\tilde{A}(m;r)$ and $\tilde{B}(r;m)$ defined by~\eqref{3-16dd}.

\begin{Corollary}\label{cb2}
With the same notation and assumptions as in Lemma~{\rm\ref{lb1}},
we have
\begin{align*}
 \big|\tilde{A}(m;r)\big|
 \begin{cases}
 < C_+ q^{\e m+\o u}, & m\ge 0,
 \\
 < C_- b^{-m} q^{\omega u}, & m<0,
 \end{cases}
 \qquad
 \big|\tilde{B}(r;m)\big|
 \begin{cases}
 < D_+ q^{\e m+(\frac12-\e-\o)u}, & m\ge 0,
 \\
 < D_- a_1^{-m}q^{u}, & m<0.
 \end{cases}
 \end{align*}
 \end{Corollary}

 Note that
 both $\tilde{A}(m;r)$ and $\tilde{B}(r;m)$ are bounded exponentially for $m$ and $r$.

\begin{proof}[Proof of Lemma~\ref{lb1}]
In this proof, $C_i$, $i=1,\dots,4$ and $D_j$, $j=1,2,3$ will denote constants which do not depend on $m$ or $r$ (i.e., $k$, $u$).

(i) First, we consider the case $m\ge 0$.
We note that the singularities of the integrand of~\eqref{3-15d} are~$0$ and $b_iq^j,$ $i=1,\dots,N$, $j\in\Z_{\ge0}$ and do not include $\tilde{a}_r^{-1}$
because of the presence of the factor $\prod_{i=1}^N(a_iz;q)_{\i}$. Thus we can extend the contour $C$ arbitrarily without changing
the value of $A(m;r)$. First we give an estimate setting the radius to be $q^{-n}$ with $n\in\Z_{>0}$ and then specify~$n$ giving an optimal estimate. With this choice of the radius $|z|=q^{-n}$, $f(m)$ in \eqref{p251-1} and some factors in the integrand in $A(m;r)$
can be estimated as
\begin{gather}
 |f(m)|<q^m,
\qquad
 \frac{2\pi q^{-n}}{|z^{1+m}|}=2\pi q^{n m},\nonumber
\\
 \left|\prod_{j=1}^N\frac{1}{(b_j/z;q)_{\i}}\right|
 <\prod_{j=1}^N\frac{1}{(b_j q^n;q)_{\i}}
 <\prod_{j=1}^N\frac{1}{(b_j;q)_{\i}},\label{bb2-3}
\end{gather}
where the factor $2\pi q^{-n}$ represents
the circumference with radius $q^{-n}$.
We rewrite the remaining factor as
\begin{align}
\label{bb2-4}
 \frac{1}{z-\tilde{a}_{r}^{-1}}
 \prod_{i=1}^N(a_i z;q)_{\i}
 =
 -a_k q^{u}
 \prod_{\substack{i=1\\i\neq k}}^N(a_i z;q)_{\i}
 \cdot (a_k z;q)_{u} \big(a_k q^{u+1}z;q\big)_{\i}.
\end{align}
The factor \smash{$-a_kq^u\prod_{\substack{i=1\\i\neq k}}^N(a_i z;q)_{\i}$} in \eqref{bb2-4} can be estimated as
\begin{gather}
 \Biggl|-a_kq^u\prod_{\substack{i=1\\i\neq k}}^N(a_i z;q)_{\i}\Biggr|\nonumber\\
\qquad < a_kq^u\prod_{\substack{i=1\\i\neq k}}^N\bigl(-a_i/q^n ;q\bigr)_{\i}
 =a_k
 q^{u-(N-1)\frac{n(n+1)}{2}}
 \prod_{\substack{i=1\\i\neq k}}^N a_i^n
 (-q/a_i;q)_n(-a_i;q)_\i
 \nonumber\\
\qquad <a_k
 q^{u-(N-1)\frac{n(n+1)}{2}}
 \prod_{\substack{i=1\\i\neq k}}^N a_i^n
 (-q/a_i;q)_\infty (-a_i;q)_\i
 <C_1 q^{u-(N-1)\frac{n(n+1)}{2}} a_{1}^{(N-1)n},\label{bb2-5}
\end{gather}
 where
 in the second equality we used the relation
 \begin{align}
 \label{bb2-6}
 \big(x/q^n;q\big)_\i=(-x)^n q^{-\frac{n(n+1)}{2}}
 (q/x;q)_{n}(x;q)_{\i}.
 \end{align}
 For the remaining factors in \eqref{bb2-4}, we find
 \begin{gather*}
 \big|(a_k z;q)_u \big(a_k q^{u+1}z;q\big)_{\i}\big| \\
 \qquad <\bigl(-a_k /q^n;q\bigr)_u\bigl(-a_k q^{u-n+1};q\bigr)_{\i}=\frac{\bigl(-a_k q^{-n};q\bigr)_\infty}{1+a_k q^{u-n}} =
 a_k^n q^{-\frac{n(n+1)}{2}}\!\frac{(-q/z_k;q)_n(-a_k;q)_{\infty}}{1+a_kq^{u-n}}\\
 \qquad <
 \begin{cases}
 C_2 a_k^n q^{-\frac{n(n+1)}{2}}, & u\ge n,
 \\
 C_3 a_k^n q^{-\frac{n(n+1)}{2}+n-u}, & n>u.
 \end{cases}
 \end{gather*}
 Here in the third equality we used \eqref{bb2-6}.
 Thus we have a bound regardless of the order of $n$ and~$u$
 \begin{align}
 \label{bb2-8}
 \big|(a_k z;q)_u\big(a_kq^{u+1}z;q\big)_{\i}\big|
 <
 C_4 a_{1}^nq^{-\frac{n(n+1)}{2}}.
 \end{align}
 Putting together \eqref{bb2-3}, \eqref{bb2-5}, and \eqref{bb2-8}, we have
 \begin{align}
 \label{bb2-9}
 |A(m;r)|&<C_+ a_{1}^{Nn} q^{nm-N\frac{n(n+1)}{2}+u+m}
 \nonumber\\
 &=C_+ a_{1}^{Nn}q^{-\frac{N}{2}\left(n-\frac{m}{N}+\frac12\right)^2+\frac{N}{2}\left(\frac{m}{N}-\frac12\right)^2+u+m}.
 \end{align}

 Now we specify $n$ which gives an optimal estimate.
 From the last expression in \eqref{bb2-9}, a~simple computation shows that such optimal bound is given by the nonnegative integer which is closest to $m/N-1/2$ and this proves the first inequality in \eqref{bb2-1}.

 Next we consider the case $m<0$. In this case, we set the radius of contour $C$ to be $b$ appearing in the assumption of (i). We find{\samepage
 \begin{align*}
 &|f(m)|<C_4
 \qquad
 \frac{2\pi b}{|z|^{1+m}}=2\pi b^{-m},
 \qquad
 \frac{1}{\big|z-\tilde{a}^{-1}_r\big|}
 <\frac{a_k q^{u}}{1-a_k bq^{u}}
 <\frac{a_k}{1-a_k} q^{\ell},
 \nonumber\\
 &\Bigg|\prod_{j=1}^M\frac{(a_iz;q)_{\i}}{(b_j/z;q)_{\i}}\Bigg|
 <\prod_{i=1}^N\frac{(-a_ib;q)_\i}{(b_i/b;q)_{\i}}.
 \end{align*}
 Combining these estimates, we arrive at the bottom estimate of \eqref{bb2-1}.}

 (ii)
 Recalling the notation $\tilde{a}_r=a_k q^{u}$ \eqref{s32-1} and calculating explicitly the residue in \eqref{3-16d}, we have
 \begin{align*}
 B(r;m)=&-a_k^{-(m+1)}q^{-u(m+1)}
 \prod_{\substack{i=1\\ i\neq k}}
 \frac{1}{(a_iq^{-u}/a_k ;q)_{\i}}\cdot
 \frac{1}{(1/q^u;q)_u (q;q)_{\i}}\cdot
 \prod_{j=1}^N(a_kb_jq^{u};q)_{\i}
 \\
 =&
 -a_k^{-(m+1)}q^{-u (m+1)+N\frac{u(u+1)}{2}}
 \prod_{\substack{i=1\\ i\neq k}}
 \frac{(-a_k/a_i)^u}{(q a_k/a_i;q)_{\i}(a_i/a_k;q)_\infty}\cdot
 \frac{(-1)^u}{(q;q)_u (q;q)_{\i}}
 \\
 &\times\prod_{j=1}^N(a_k b_jq^{u};q)_{\i},
 \end{align*}
 where in the second equality we used \eqref{bb2-6}.
 Using the condition $a_1/a_N< q^{-1/2+\e}$, we have the following estimate,
 \begin{align}
 \label{bb2-12}
 |B(r;m)|&<D_1 a_k^{-(m+1)} q^{-u m +N\frac{u(u+1)}{2}-N\left(\frac12-\e\right)u -\left(\frac12+\e\right)u}.
 \end{align}
 For $m\ge 0$, combining \eqref{bb2-12} with
 the following bound
 \begin{align*}
 q^{-u m+N\frac{u(u+1)}{2}-N\left(\frac12-\e\right)u}
 =
 q^{\frac{N}{2}\left(u-\frac{m}{N}+\e\right)^2-\frac{N}{2}\left(\frac{m}{N}-\e\right)^2}
 <D_2
 q^{-\frac{m^2}{2N}+\e m},
 \end{align*}
 we have the desired estimate. On the other hand
 for $m<0$, we use
 \begin{align}
 \label{bb2-14}
 q^{-u m}<q^{u},
 \qquad
 q^{N\frac{u(u+1)}{2}-N\left(\frac12-\e\right)u -\left(\frac12+\e\right)u}<D_3.
 \end{align}
 From \eqref{bb2-12} and \eqref{bb2-14}, we obtain
 the bottom relation in \eqref{bb2-2}.
\end{proof}

\section[Decay estimates of q-Pochhammer symbols]{Decay estimates of $\boldsymbol{q}$-Pochhammer symbols}
\label{ac}
In this section, we provide simple bounds for
the $q$-Pochhammer symbol $(z;q)_{\i}$, where
$z\in\C$ and $0<q<1$. For our purpose, it is sufficient
to have bounds in the following three cases though
they are not exclusive nor exhaustive.
\begin{Lemma}\label{lc1}
 Let us write $z=a+\sqrt{-1}b$ with $a,b\in\R$.
 For $0<q<1$ fixed, there exist positive constants
 $c_1=c_1(q)$ and $c_2=c_2(q)$ such that
 \begin{itemize}\itemsep=0pt
 \item[\rm{(i)}] for $a\le -1$,
 \begin{align}
 \label{lc1-1}
 |(z;q)_\i|\ge c_1 \exp \big[c_2\log^2 |a|\big],
 \end{align}
 \item[\rm{(ii)}] for $a<1$ and $|b|>1$,
 \begin{align}
 \label{lc1-2}
 |(z;q)_\i|\ge c_1 \exp\big[c_2\log^2 |b|\big],
 \end{align}
 \item[\rm{(iii)}] for $a>1$,
 the following two types of lower bounds
 hold
 \begin{align}
 \label{lc1-3}
 |(z;q)_\i|
 &\ge c_1 b^2 a^{-\frac32}
 \exp\big[c_2\log^2 a\big],\\
 \label{lc1-3'}
 |(z;q)_\i|
 &\ge c_1 \big(q^{\a-1}-1\big)\big(1-q^{\a}\big)
 \exp\big[c_2\log^2 a\big],
 \end{align}
 where $\a=\left\lceil \log_{q^{-1}}a\right\rceil-\log_{q^{-1}}a$
 with $\lceil x\rceil$ being the minimum integer
 greater than $x$.
 \end{itemize}
\end{Lemma}

\begin{proof}
 We repeatedly use the simple estimate
 \begin{align}
\label{lc1-5}
 |1-w|=\sqrt{(1-\Re w)^2+(\Im w)^2}\ge
 \begin{cases}
 |1-\Re w|,\\
 |\Im w|,
 \end{cases}
\end{align}
 which hold for any $w\in\C$. We use the top or bottom estimate depending on the situations.

 (i) In terms of $\a\in [0,1)$ and $J\in\Z_{>0}$, we express $a(<-1)$
 as $a=-q^{\a-J}$.
 Applying the top estimate in \eqref{lc1-5}
 to the factor $\big|1-z q^k\big|$, we have
 \begin{align}
 \label{lc1-6}
 \big|1-z q^k\big|=\sqrt{\big(1+aq^{k}\big)^2+\big(bq^k\big)^2}
 \ge 1+q^{\a-J+k}
 \ge
 \begin{cases}
 q^{\a-J+k},
 \\
 1-q^{k+1-J},
 \end{cases}
 \end{align}
 where both estimates holds for $\forall k=0,1,2,\dots$. Choosing the top (resp.\ bottom) estimate for~${k<J}$ (resp.\ $k\ge J$), we have
 \begin{align}
 \label{lc1-7}
 |(z;q)_\i|\ge q^{\a J-\frac{J(J+1)}{2}}(q;q)_{\i}
 =
 q^{-\frac{1}{2}(\a-J)^2+\frac{1}{2}(\a-J)
 +\frac{\a^2-\a}{2}}(q;q)_{\i}
 \ge
 q^{-\frac{1}{2}(\a-J)^2}(q;q)_{\i},
 \end{align}
 where in the last inequality we used the fact
 $0\le\a<1$ and $1\le J$. Recalling $\a-J=\log_q|a|$,
 one easily sees
 \begin{align}
 \label{lc1-8}
 q^{-\frac{(\a-J)^2}{2}}=
 {\rm e}^{\frac{-(\log |a|)^2}{2\log q}}.
 \end{align}
 From \eqref{lc1-7} and \eqref{lc1-8},
 we obtain \eqref{lc1-1}.

 (ii)
 We adopt the same strategy as above.
 In this case we write $|b|=q^{\a-J}$, where
 $0\le \a<1$ and $J\in\Z_{>0}$. Using \eqref{lc1-5},
 we have for $k=0,1,2,\dots$
 \begin{align}
 \label{lc1-9}
 \big|1-z q^k\big|
 \ge
 \begin{cases}
 |b|q^{k}=q^{\a-J+k},
 \\
 1-aq^k\ge 1-q^{k+1-J}.
 \end{cases}
 \end{align}
 Applying the top (resp.\ bottom) estimate in \eqref{lc1-9} to the case $k<J$ (resp.\ $k\ge J$) and noting that the estimate \eqref{lc1-9}
 is exactly the same as \eqref{lc1-6} except
 that $\a-J$ represents $\log_q|b|$, we get~\eqref{lc1-2} similarly to the case
 (i).

 (iii)
 As in the case (i),
 we express $a$ as $a=q^{\a-J}$ with $a\in[0,1)$ and $J\in\{1,2,\dots\}$.
 First, we prove \eqref{lc1-3}.
 We choose the estimates as follows:
 \begin{gather}
 \big|1-z q^k\big|\nonumber
 \\
 \qquad\ge
 \begin{cases}
 \big|1-a q^k\big|=q^{\a-J+k}\big(1-q^{J-\a-k}\big)&\\
 \phantom{ \big|1-a q^k\big|}{} \ge q^{\a-J+k}\big(1-q^{J-k-1}\big), & {\text{for $0\le k\le J-2$}},
 \\
 bq^k, & {\text{for $k=J-1,J$}},
 \\
 \big|1-a q^k\big|=1-q^{\a-J+k}\ge 1-q^{-J+k}, & {\text{for $J+1\le k$}}.
 \end{cases}\label{lc1-10}
 \end{gather}
 Thus we have
 \begin{align*}
 |(z;q)_{\i}|
 &\ge
 \prod_{k=0}^{J-2} q^{\a-J+k}\big(1-q^{J-k-1}\big)
 \cdot
 \prod_{k=J-1}^{J} bq^{k}
 \cdot
 \prod_{k=J+1}^{\i} \big(1-q^{-J+k}\big)
 \\
 &=
 q^{\a(J-1)-\frac{J(J+1)}{2}+2J}b^2(q;q)_{J-1}(q;q)_{\infty}\ge
 q^{-\frac{(J-\alpha)^2}{2}+\frac32(J-\alpha)+\frac{\alpha^2}{4}+\frac{3\alpha}{2}}b^2(q;q)_{\infty}^2.
 \end{align*}
 Thus from $a=q^{\alpha-J}$ and \eqref{lc1-8}, we
 obtain the desired form \eqref{lc1-3}.

 The other estimate \eqref{lc1-3'}
 can be readily
 obtained by replacing the second inequality in \eqref{lc1-10}
 (corresponding to the case $k=J-1, J$)
 with $\big|1-zq^k\big|\ge \big|1-aq^k\big|=\big|1-q^{\a-J+k}\big|$.
\end{proof}

\subsection*{Acknowledgements}
 This work was initiated
 during the MATRIX program "Non-equilibrium systems and special functions" in 2018. The authors are grateful to Dan Betea and J{\'e}r{\'e}mie Bouttier for
 introduction to their recent work \cite{BeteaBouttier2019PeriodicSchur}
 and for
 discussions about the periodic Schur measure and the free Fermion at positive
 temperature during the program.
 In particular, TS thanks DB for pointing out~I.5,28(a) of Macdonald's book. We greatly appreciate creative atmosphere and warm hospitality of the MATRIX institute. We are grateful to Alexei Borodin and Guillaume Barraquand whose conversations with led to Section~\ref{sec:comparison} and to a better understanding of relations between our result and those of \cite{Borodin2018momentmatch}.
 We thank the anonymous referees for careful reading of the manuscript and helpful comments.
 The work of TI has been supported by JSPS KAKENHI Grant No.\ JP16K05192, No.\ JP19H01793, No.\ JP20K03626, and No.\ JP22H01143.
 The work of TS has been supported by JSPS KAKENHI Grants No.\ JP15K05203, No.\ JP16H06338, No.\ JP18H01141, No.\ JP18H03672, No.\ JP19L03665, No.\ JP21H04432, No.\ JP22H01143.
 The work of MM has been partially supported by the European Union’s Horizon 2020 research and innovation programme under the Marie Sklodowska-Curie grant agreement No. 101030938.

\pdfbookmark[1]{References}{ref}
\LastPageEnding

\end{document}